\documentclass[lettersize,journal]{IEEEtran}
\usepackage{amsmath,amsfonts,amssymb,amsthm}

\usepackage{algorithmic}
\usepackage{algorithm}
\usepackage{array}
\usepackage{textcomp}
\usepackage{url}
\usepackage{verbatim}
\usepackage{cite}
\usepackage{booktabs}

\usepackage{tikz}
\usepackage{tikz-3dplot}
\usepackage{xcolor}

\usepackage[font=footnotesize,labelfont=rm,textfont=rm]{subfig}

\usepackage{bm}
\usepackage{bbm}
\usepackage{upgreek}
\usepackage{multirow}

\usepackage{psfrag}
\usepackage{epsfig}
\usepackage{epstopdf}

\usepackage{mathrsfs}
\usepackage{courier}
\usepackage{tikz}

\usepackage{graphicx}
\usepackage{listings}
\usepackage{xcolor} 

\usepackage[font=footnotesize]{caption}
\DeclareCaptionLabelFormat{none}{}
\usetikzlibrary{positioning}

\newtheorem{theorem}{Theorem}
\newtheorem{lemma}{Lemma} 
\newtheorem{corollary}{Corollary} 
\newtheorem{remark}{Remark}%
\newtheorem{definition}{Definition}

\DeclareCaptionLabelFormat{algnonumber}{Algorithm}
\makeatletter
\newsavebox\myboxA
\newsavebox\myboxB
\newlength\mylenA
\newcommand*\xoverline[2][0.75]{%
	\sbox{\myboxA}{$\m@th#2$}%
	\setbox\myboxB\null
	\ht\myboxB=\ht\myboxA%
	\dp\myboxB=\dp\myboxA%
	\wd\myboxB=#1\wd\myboxA
	\sbox\myboxB{$\m@th\overline{\copy\myboxB}$}
	\setlength\mylenA{\the\wd\myboxA}
	\addtolength\mylenA{-\the\wd\myboxB}%
	\ifdim\wd\myboxB<\wd\myboxA%
	\rlap{\hskip 0.5\mylenA\usebox\myboxB}{\usebox\myboxA}%
	\else
	\hskip -0.5\mylenA\rlap{\usebox\myboxA}{\hskip 0.5\mylenA\usebox\myboxB}%
	\fi}
\makeatother

\usepackage{xcolor}
\usepackage{xpatch}
\makeatletter
\ExplSyntaxOn
\cs_new:Npn \bibColoredItems #1#2
{
	\clist_map_inline:nn {#2} { \cs_new:cpn {bib@colored@##1} {#1} }
}
\ExplSyntaxOff

\newcommand\bib@setcolor[1]{%
	\ifcsname bib@colored@#1\endcsname
	\expandafter\color\expandafter{\csname bib@colored@#1\endcsname}
	\else
	\normalcolor
	\fi
}

\xpatchcmd\@bibitem
{\item}
{\bib@setcolor{#1}\item}
{}{\fail}

\xpatchcmd\@lbibitem
{\item}
{\bib@setcolor{#2}\item}
{}{\fail}
\makeatother

\begin{document}

\title{Revisiting time-variant complex conjugate matrix equations with their corresponding real field time-variant large-scale linear equations, neural hypercomplex numbers space compressive approximation approach}

\author{
Jiakuang He, Dongqing Wu

\thanks{This work is aided by the Project Supported by the Guangzhou Science and Technology Program (with number 2023E04J1240).
	(Corresponding author: Dongqing Wu)}
\thanks{Jiakuang He is graduated from Guangzhou University of Chinese
	Medicine, Guangzhou 510006, P. R. China; Jiakuang He and Dongqing Wu are with School of Mathematics and Data Science, Zhongkai University of Agriculture and Engineering, Guangzhou 510225, P.
	R. China (e-mails: rickwu@zhku.edu.cn).}}



\maketitle

\begin{abstract}
Large-scale linear equations and high dimension have been hot topics in deep learning, machine learning, control, and scientific computing. Because of special conjugate operation characteristics, time-variant complex conjugate matrix equations need to be transformed into corresponding real field time-variant large-scale linear equations. In this paper, zeroing neural dynamic models based on complex field error (called Con-CZND1) and based on real field error (called Con-CZND2) are proposed for in-depth analysis. Con-CZND1 has fewer elements because of the direct processing of complex matrices. Con-CZND2 needs to be transformed into the real field, with more elements, and its performance is affected by the main diagonal dominance of coefficient matrices. A neural hypercomplex numbers space compressive approximation approach (NHNSCAA) is innovatively proposed. Then Con-CZND1\_conj model is constructed. Numerical experiments verify Con-CZND1\_conj model effectiveness and highlight NHNSCAA importance.

\end{abstract}


\section{Introduction}
\begin{figure*}[!t]
	\centering
	\resizebox{0.6\textwidth}{!}{%
	\begin{tikzpicture}[
		every node/.style={
			draw, 
			circle, 
			align=center, 
			inner sep=2pt, 
			minimum size=1.5cm 
		},
		edge from parent/.style={
			->,
			thick,
			draw,
			shorten >=1pt
		},
		myarrow/.style={
			->,
			>=latex,
			shorten >=1pt,
			shorten <=1pt,
			draw
		},
		level 1/.style={sibling distance=6cm},
		level 2/.style={sibling distance=5cm},
		]
		
		\node[draw, circle, align=center, inner sep=2pt, fill=blue!20] (center) {Large-scale\\linear\\equations};
		
		\node[draw, circle, align=center, inner sep=2pt, fill=green!20] (iter) [below right=of center] {\textit{Iterative}\\\textit{methods}};
		\node[draw, circle, align=center, inner sep=2pt, fill=yellow!20] (decomp) [below left=of center] {\textit{Direct}\\\textit{decomposition}\\\textit{methods}};
		\node[draw, circle, align=center, inner sep=2pt, fill=orange!20] (sparse) [left=of center] {\textit{Sparse}\\\textit{matrix}\\\textit{techniques}};
		\node[draw, circle, align=center, inner sep=2pt, fill=red!20] (preproc) [above left=of center] {\textit{Preprocessing}\\\textit{techniques}\\\textit{\cite{6016756}}};
		\node[draw, circle, align=center, inner sep=2pt, fill=purple!20] (parallel) [above right=of center] {\textit{Parallel}\\\textit{and}\\\textit{distributed}\\\textit{computations}\\\textit{\cite{6926792}}};
		
		\draw[myarrow] (center) -- (iter);
		\draw[myarrow] (center) -- (decomp);
		\draw[myarrow] (center) -- (sparse);
		\draw[myarrow] (center) -- (preproc);
		\draw[myarrow] (center) -- (parallel);
		
		\draw[myarrow] (iter) -- ++(-1,-1) node[left] {\textit{Jacobi}\\\textit{iteration}\\\textit{\cite{8350317}}};
		\draw[myarrow] (iter) -- ++(1,1) node[above right] {\textit{Gauss-Seidel}\\\textit{iteration}\\\textit{\cite{9325945}}};
		\draw[myarrow] (iter) -- ++(1,-1) node[below right] {\textit{Successive}\\\textit{over}\\\textit{relaxation}\\\textit{(SOR)}\\\textit{\cite{9426571}}};
		
		\draw[myarrow] (decomp) -- ++(-1,-1.5) node[below left] {\textit{LU}\\\textit{factorization}\\\textit{\cite{10143289}}};
		
		\draw[myarrow] (sparse) -- ++(-1.5,1) node[above left] {\textit{Compressing}\\\textit{row}\\\textit{and}\\\textit{column}\\\textit{storage}\\\textit{\cite{6574845}}};
		
	\end{tikzpicture}}
	\caption{Recent main solving methods of LSLE.}
	\label{fig.mainlslesolving}
\end{figure*}
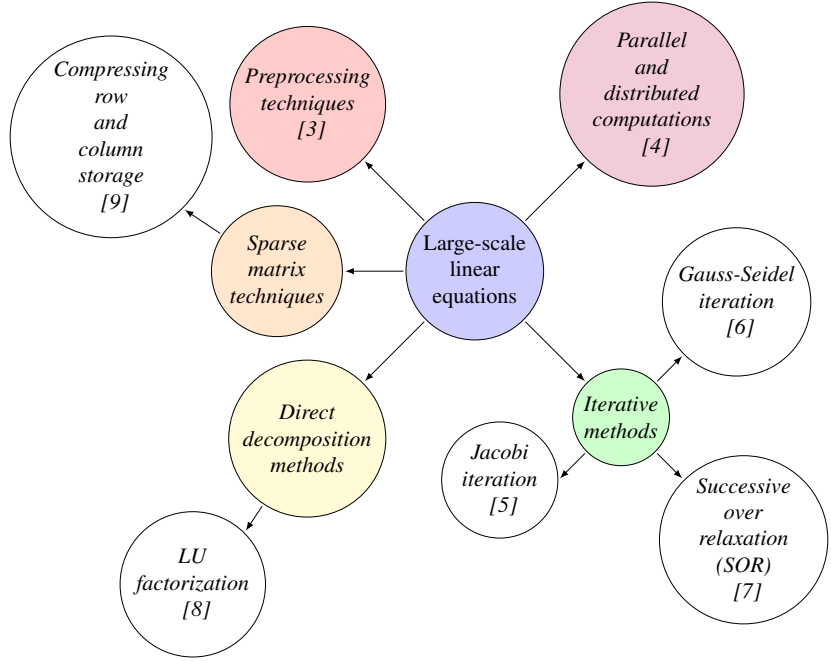

\IEEEPARstart{L}{arge}-scale linear equations (LSLE)
\cite{167018}
are a generalization of 
linear equations
\cite{WU2023435}.
Because of the large computation, it has become a hot spot of the current study. Recently, the main solving methods of LSLE are shown in Fig. \ref{fig.mainlslesolving}.
The corresponding real field time-variant LSLE (RFTVLSLE) becomes more difficult to solve because of time-variant elements. According to known studies, there is no well-suited entry point for optimization, which is a bottleneck for deep learning
\cite{LIAO2022440},
machine learning
\cite{8730514},
control
\cite{10127863},
and scientific computing
\cite{9585655}.
An early time-variant version of complex conjugate matrix equation (CCME) 
\cite{Wu2017}
is time-variant standard Sylvester-conjugate matrix equations (TVSSCME)
\cite{He2024ZeroingND},
which is also being studied in depth. Because of the conjugate computation, the solved complex matrices have to be separated by real transformations. Con-CZND1
\cite{He2024ZeroingND}
model with based complex field error and Con-CZND2
\cite{He2024ZeroingND}
model with based real field error are also proposed simultaneously.
All of the above models have to be transformed into RFTVLSLE. But ultimately, Con-CZND2 is affected by this, while Con-CZND1 is not.

High dimension (HD)
\cite{Hastie2009}
is the essential reason. HD not only causes an increase in computational complexity but also brings ``curse of dimensionality” (COD)
\cite{Hastie2009},
resulting in large errors. Therefore, some dimensional reduction techniques must be developed urgently. And in the optimization of HD matrix computation, the time-variant problem must also be taken into account. However, according to the known study, it is still in the direction of time-invariant extension. Considering the real world, the time-variant study must also pay attention to.

Zeroing Neural Networks (ZNN)
\cite{1031938,uhligZhangNeuralNetworks2024}
have been proposed to solve TVSSCME as a class of neural networks (NN) based on ordinary differential equations (ODE). It cleverly exploits the solution of function matrices through parallel sub-element descent, avoiding the difficulty of gradient derivation and real representation \cite{4203075}.
However, this zeroing neural dynamics (ZND)
\cite{He2024ZeroingND}
is currently in a preliminary stage for the time-variant version of CCME and relevant knowledge backgrounds are not well developed. Moreover, a systematic theory is urgently needed.

Errors in NN mainly come from activate function errors,
sampling discretion errors,
and space compressive approximation errors. The activation function under ZND is currently developed to be noise-resistant
\cite{article},
finite time
\cite{2024fix},
and adjustable factor
\cite{9078848}.
Sampling discretion errors are currently developed to more than 10 points sampling 
\cite{WU2022391,10355921,WU202344}.
The study of space compressive approximation, on the other hand, is the most difficult and involves the nature of hypercomplex numbers (HN)
\cite{Valle2024UniversalAT}.
As there is not a specific representative time-variant problem. Therefore, in this paper, a systematic study of space compressive approximation is carried out in conjunction with TVSSCME,
which is being studied. From a new perspective, a new dimensional descent approach is proposed for optimizing NN. The corresponding RFTVLSLE and HD problems are handled simultaneously. It not only complements the perspective of space compressive approximation errors, but also provides a detailed description of the universal approximation theorem of NN and refines the corresponding theory.
 
The rest of the paper is organized as follows: Section \uppercase\expandafter{\romannumeral2} takes an novel viewpoint of directly using Con-CZND2 model to solve the insignificant main diagonal dominance RFTVLSLE. This RFTVLSLE is equivalent to TVSSCME, and the necessary background for this paper is also given. Section \uppercase\expandafter{\romannumeral3} explains the deduction of solving RFTVLSLE, proposes NHNSCAA, and constructs the corresponding conjugate model Con-CZND1\_conj for Con-CZND1. Section \uppercase\expandafter{\romannumeral4} designs simulations and special experiments to verify the validity of Con-CZND1 model and Con-CZND1\_conj model. The differences between those models are also compared to highlight the effectiveness of NHNSCAA. Finally, a conclusion and future directions are given with appendix references. Before starting the next section, the main contributions of this paper are as follows.
\begin{enumerate}
	\item Con-CZND2 model is explained from a novel viewpoint. Then, necessary background knowledge is provided to highlight the similarities and differences between RFTVLSLE and TVSSCME.
	
	\item With the idea of space compressive approximation, a new dimensional descent approach is proposed. On this basis, the conjugate model Con-CZND1\_conj corresponding to Con-CZND1 model is constructed. Then, two models can solve RFTVLSLE corresponding to TVSSCME.
	 
	\item Con-CZND1 model and Con-CZND1\_conj model can be used to verify he significance of space compressive approximation in neural networks. These models highlight the details of the universal approximation theorem.
\end{enumerate} 

\section{Problem Formulation}
Consider a RFTVLSLE \eqref{eq.RFTVLSLE}
\cite{bevisConsimilarityMatrixEquation1987,wuSolutionsMatrixEquations2006,Wu2017,He2024ZeroingND}
as
\begin{equation} \label{eq.RFTVLSLE}
	\begin{split}
		\begin{bmatrix}
			K_{11}(\tau)	&K_{12}(\tau) \\
			K_{21}(\tau)	&K_{22}(\tau)
		\end{bmatrix}
		\begin{bmatrix}
			\mathrm{vec}(X_{\mathrm{r}}(\tau))\\
			\mathrm{vec}(X_{\mathrm{i}}(\tau))	
		\end{bmatrix}
		=\begin{bmatrix}
			\mathrm{vec}(C_{\mathrm{r}}(\tau))\\
			\mathrm{vec}(C_{\mathrm{i}}(\tau))	
		\end{bmatrix},
	\end{split}
\end{equation}
where $\tau \ge 0$ represents the real-time, $W_{\mathrm{R}}(\tau)=[K_{11}(\tau),K_{12}(\tau);
K_{21}(\tau),K_{22}(\tau)]\in
\mathbb{R}^{2mn\times 2mn}$, $K_{11}(\tau)=(F_{\mathrm{r}}^{\mathrm{T}}(\tau) \otimes I_{m})-(I_{n} \otimes A_{\mathrm{r}}(\tau))\in
\mathbb{R}^{nm\times mn}$, $K_{12}(\tau)=-(F_{\mathrm{i}}^{\mathrm{T}}(\tau) \otimes I_{m}+I_{n} \otimes A_{\mathrm{i}}(\tau))\in
\mathbb{R}^{nm\times mn}$,
$K_{21}(\tau)=(F_{\mathrm{i}}^{\mathrm{T}}(\tau) \otimes I_{m})-(I_{n} \otimes A_{\mathrm{i}}(\tau))\in
\mathbb{R}^{nm\times mn}$,
$K_{22}(\tau)=(F_{\mathrm{r}}^{\mathrm{T}}(\tau) \otimes I_{m}+I_{n}\otimes A_{\mathrm{r}}(\tau))\in
\mathbb{R}^{nm\times mn}$,
$X_{\mathrm{R}}(\tau)=[\mathrm{vec}(X_{\mathrm{r}}(\tau));
\mathrm{vec}(X_{\mathrm{i}}(\tau))]\in\mathbb{R}^{2mn\times 1}$ is a time-variant vector to be computed,
$B_{\mathrm{R}}(\tau)=[\mathrm{vec}(C_{\mathrm{r}}(\tau)); \mathrm{vec}(C_{\mathrm{i}}(\tau))]\in\mathbb{R}^{2mn\times 1}$.
For simplicity, RFTVLSLE \eqref{eq.RFTVLSLE} is written as \eqref{eq.simplify.RFTVLSLE}:
\begin{equation} \label{eq.simplify.RFTVLSLE}
	W_{\mathrm{R}}(\tau)X_{\mathrm{R}}(\tau) =B_{\mathrm{R}}(\tau).
\end{equation}

And RFTVLSLE \eqref{eq.RFTVLSLE} corresponding TVSSCME \eqref{eq.TVSSCME}
\cite{bevisConsimilarityMatrixEquation1987,wuSolutionsMatrixEquations2006,Wu2017,He2024ZeroingND}
is
\begin{equation} \label{eq.TVSSCME}
	\begin{split}
		X(\tau)F(\tau)-A(\tau)\overline{X}(\tau)=C(\tau),
	\end{split}
\end{equation}
where $F(\tau)\in
\mathbb{C}^{n\times n}$, $A(\tau)\in
\mathbb{C}^{m\times m}$, and $C(\tau)\in
\mathbb{C}^{m\times n}$ are known as time-variant matrices,
$X(\tau)\in
\mathbb{C}^{m\times n}$ is a time-variant matrix to be computed.

Before deduction, some relevant background knowledge and formulas are added:
\begin{definition}
	Where $\tau \ge 0$ represents the real-time, time-variant complex matrices elements
	\cite{He2024ZeroingND}
	are defined as follows:
\begin{equation}\label{eq.define.complexmatrix.elements}
	\tilde{m}_{st}(\tau)=m_{\mathrm{r},st}(\tau)+\mathrm{i}m_{\mathrm{i},st}(\tau),
\end{equation}
where $s\in
\mathbb{I}[1,p]$, $t\in
\mathbb{I}[1,q]$, $\mathbb{I}[m,n]$ means the set of integers from m to n, $m_{\mathrm{r},st}(\tau)\in
\mathbb{R}$ is a real coefficient, $m_{\mathrm{i},st}(\tau)\in
\mathbb{R}$ is a imaginary coefficient, $\mathrm{i}$ is an imaginary unit, same as below. Also, if $m_{\mathrm{i},st}(\tau)=0$, the elements $\tilde{m}_{st}(\tau)=\left (m_{\mathrm{r},st}(\tau)+\mathrm{i}\times0\right )\in
\mathbb{R}$. For the sake of uniformity, the elements are represented by $m_{st}(\tau)=\tilde{m}_{st}(\tau)$ under the real matrices.

Correspondingly, the time-variant complex matrices
\cite{He2024ZeroingND}
can be written in the following form
\begin{equation}\label{eq.define.complexmatrix}
	M(\tau)=M_{\mathrm{r}}(\tau)+\mathrm{i}M_{\mathrm{i}}(\tau),
\end{equation}	
where $M(\tau)\in
\mathbb{C}^{p\times q}$ is any complex matrix, $M_{\mathrm{r}}(\tau)\in
\mathbb{R}^{p\times q}$ is the real coefficient matrix of $M(\tau)$, $M_{\mathrm{i}}(\tau)\in
\mathbb{R}^{p\times q}$ is the imaginary coefficient matrix of $M(\tau)$.
In Fig. \ref{fig.complexmatrix},
The structure of \eqref{eq.define.complexmatrix} can be seen schematically.
The conjugate matrix corresponding to $M(\tau)$ is $\overline{M}(\tau)=M_{\mathrm{r}}(\tau)-\mathrm{i}M_{\mathrm{i}}(\tau)$, where $\overline{M}(\tau)\in\mathbb{C}^{p\times q}$.
\end{definition}
	\begin{figure}[!t]
	\centering	
	\begin{tikzpicture}  
		\tikzstyle{grid}=[gray,thin]
		\begin{scope}[xshift=1.75cm, yshift=2.25cm]  
			\foreach \x in {0,1.5,3} {  
				\foreach \y in {0,1.5,3} {  
					\draw[grid,dashed] (\x,\y) rectangle (\x+1.5,\y+1.5);  
				}  
			}  
			
			\node at (0.75,0.75) {$\tilde{m}_{p1}(\tau)$};  
			\node at (2.25,0.75) {$\cdots$};  
			\node at (3.75,0.75) {$\tilde{m}_{pq}(\tau)$};  
			\node at (0.75,2.25) {$\vdots$};  
			\node at (2.25,2.25) {$\ddots$};  
			\node at (3.75,2.25) {$\vdots$};  
			\node at (0.75,3.75) {$\tilde{m}_{11}(\tau)$};  
			\node at (2.25,3.75) {$\cdots$};  
			\node at (3.75,3.75) {$\tilde{m}_{1q}(\tau)$};

		\end{scope}
		\draw[->, thick] (4,2) -- (4,1);       
		\begin{scope}[xshift=4cm, yshift=-3.75cm]  
			\foreach \x in {0,1.5,3} {  
				\foreach \y in {0,1.5,3} {  
					\draw[grid] (\x,\y) rectangle (\x+1.5,\y+1.5);  
				}  
			}  
			
			\node at (0.75,0.75) {$m_{\mathrm{i},p1}(\tau)$};  
			\node at (2.25,0.75) {$\cdots$};  
			\node at (3.75,0.75) {$m_{\mathrm{i},pq}(\tau)$};  
			\node at (0.75,2.25) {$\vdots$};  
			\node at (2.25,2.25) {$\ddots$};  
			\node at (3.75,2.25) {$\vdots$};  
			\node at (0.75,3.75) {$m_{\mathrm{i},11}(\tau)$};  
			\node at (2.25,3.75) {$\cdots$};  
			\node at (3.75,3.75) {$m_{\mathrm{i},1q}(\tau)$};

		\end{scope}
		\node[anchor=east, font=\bfseries\large] at (3.75,-1.6225) {$+$ $\mathrm{i}$}; 
		\begin{scope}[xshift=0cm, yshift=-7.75cm]  
			\foreach \x in {0,1.5,3} {  
				\foreach \y in {0,1.5,3} {  
					\draw[grid,fill=lightgray] (\x,\y) rectangle (\x+1.5,\y+1.5);  
				}  
			}  
			
			\node at (0.75,0.75) {$m_{\mathrm{r},p1}(\tau)$};  
			\node at (2.25,0.75) {$\cdots$};  
			\node at (3.75,0.75) {$m_{\mathrm{r},pq}(\tau)$};  
			\node at (0.75,2.25) {$\vdots$};  
			\node at (2.25,2.25) {$\ddots$};  
			\node at (3.75,2.25) {$\vdots$};  
			\node at (0.75,3.75) {$m_{\mathrm{r},11}(\tau)$};  
			\node at (2.25,3.75) {$\cdots$};  
			\node at (3.75,3.75) {$m_{\mathrm{r},1q}(\tau)$};

		\end{scope}
		\draw[->, thick] (4,-8) -- (4,-9);  
		
		\begin{scope}[xshift=3.5cm, yshift=-10.75cm] 
			[every node/.style={font=\small}]  
			\coordinate (A) at (0,0,0);  
			\coordinate (B) at (1.5,0,0);  
			\coordinate (C) at (1.5,1.5,0);  
			\coordinate (D) at (0,1.5,0);  
			\coordinate (A1) at (0,0,1);
			\coordinate (A11) at (0,0,0.5);    
			\coordinate (B1) at (1.5,0,1);
			\coordinate (B11) at (1.5,0,0.5); 
			\coordinate (C1) at (1.5,1.5,1);
			\coordinate (C11) at (1.5,1.5,0.5);  
			\coordinate (D1) at (0,1.5,1);
			\coordinate (D11) at (0,1.5,0.5);   
			
			\fill[gray, opacity=0.5] (A1) -- (B1) -- (C1) -- (D1) -- cycle; 
			\fill[gray, opacity=0.5] (C1) -- (C11) -- (D11) -- (D1) -- cycle; 
			\fill[gray, opacity=0.5] (C1) -- (C11) -- (B11) -- (B1) -- cycle; 
			
			\draw(A) -- node[pos=0.5,below] {} (B);  
			\draw(B) -- node[pos=0.5,right] {$\mathit{p}$} (C);  
			\draw(C) -- node[pos=0.5,right] {} (C1);  
			\draw(C1) -- node[pos=0.5,above] {} (B1);  
			\draw(B1) -- node[pos=0.5,below] {$\mathit{q}$} (A1);  
			\draw(A1) -- node[pos=0.5,left] {} (D1);  
			\draw(D1) -- node[pos=0.5,left] {} (D);  
			\draw(D) -- node[pos=0.5,below] {} (A);  
			
			\draw(A) -- node[pos=0.5,left] {} (A1);  
			\draw(B) -- node[pos=0.4,below right] {$2$} (B1);  
			\draw(C) -- node[pos=0.5,right] {} (D);  
			\draw(D1) -- node[pos=0.5,above] {} (C1);
		\end{scope}    
	\end{tikzpicture}  
	\caption{Schematic figure of \eqref{eq.define.complexmatrix} structure.}
	\label{fig.complexmatrix}  
\end{figure}
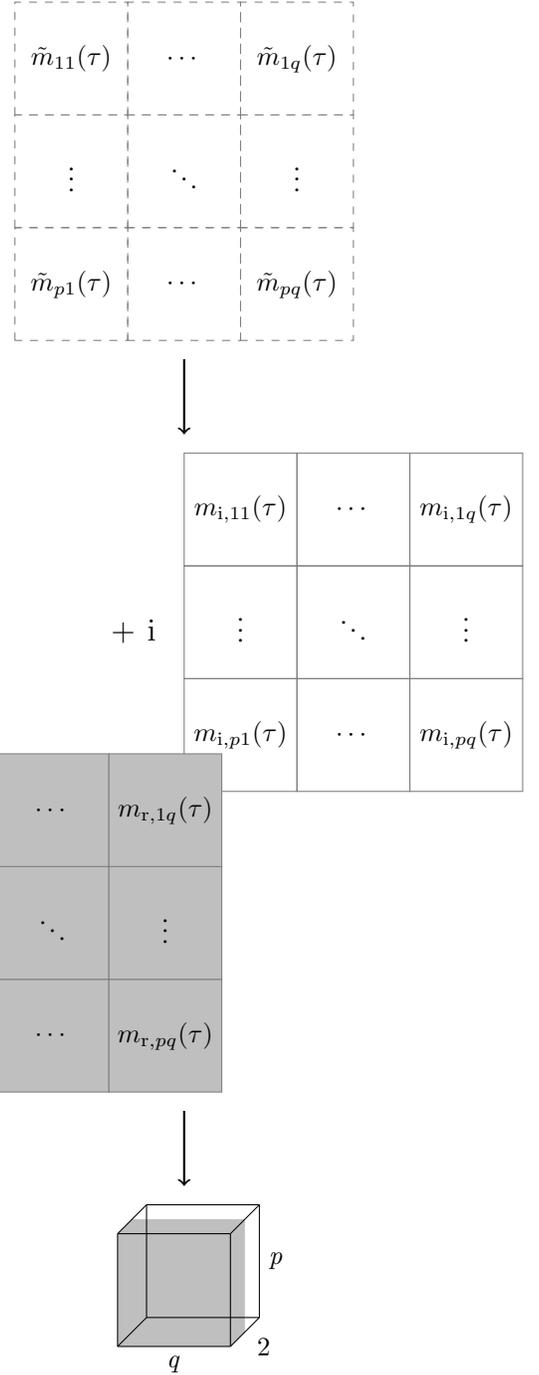 
\begin{lemma}
Where $\tau \ge 0$ represents the real-time, $A(\tau)\in
\mathbb{C}^{m\times n}$, $B(\tau)\in
\mathbb{C}^{s\times t}$, $X(\tau)\in
\mathbb{C}^{n\times s}$, are time-variant matrices,
the following equation can be obtained:
\begin{equation}\label{eq.infer.complexkroneckerproductvectorization} 
	\mathrm{vec}(A(\tau)X(\tau)B(\tau))=(\overline{B^{\mathrm{H}}}(\tau)\otimes A(\tau))\mathrm{vec}(X(\tau)).	
\end{equation}
\begin{proof}
	The proof can be generalized from Theorem 1 of \cite{He2024ZeroingND}.
\end{proof}
\begin{corollary}
	Where $A(\tau)\in
	\mathbb{R}^{m\times n}$, $B(\tau)\in
	\mathbb{R}^{s\times t}$, $X(\tau)\in
	\mathbb{R}^{n\times s}$, and $\tau \ge 0$ represents the real-time, \eqref{eq.infer.complexkroneckerproductvectorization} converts to: \begin{equation}\label{eq.infer.realkroneckerproductvectorization} 
		\mathrm{vec}(A(\tau)X(\tau)B(\tau))=(B^{\mathrm{T}}(\tau)\otimes A(\tau))\mathrm{vec}(X(\tau)).	
	\end{equation}
\end{corollary}
\end{lemma}
Based on the above background knowledge and equations, \eqref{eq.TVSSCME} derives \eqref{eq.RFTVLSLE} as follows:	
\begin{corollary}
According to \eqref{eq.define.complexmatrix}, TVSSCME \eqref{eq.TVSSCME} is firstly split to obtain the following:
\begin{equation}\label{eq.transccsme.variant}
	\begin{split}
	(X_{\mathrm{r}}(\tau)+\mathrm{i}X_{\mathrm{i}}(\tau))(F_{\mathrm{r}}(\tau)+\mathrm{i}F_{\mathrm{i}}(\tau))\\-(A_{\mathrm{r}}(\tau)+\mathrm{i}A_{\mathrm{i}}(\tau))(X_{\mathrm{r}}(\tau)-\mathrm{i}X_{\mathrm{i}}(\tau))\\=(C_{\mathrm{r}}(\tau)+\mathrm{i}C_{\mathrm{i}}(\tau)).
	\end{split}	
\end{equation}

Separating \eqref{eq.transccsme.variant} is then carried out to obtain the equivalent real field \eqref{eq.dividetransccsme.variant} as follows:
\begin{equation}\label{eq.dividetransccsme.variant}
	\begin{split}
		\left\{\begin{matrix}
		X_{\mathrm{r}}(\tau)F_{\mathrm{r}}(\tau)-X_{\mathrm{i}}(\tau)F_{\mathrm{i}}(\tau)-A_{\mathrm{r}}(\tau)X_{\mathrm{r}}(\tau)-A_{\mathrm{i}}(\tau)X_{\mathrm{i}}(\tau)\\=C_{\mathrm{r}}(\tau),	\\
		X_{\mathrm{i}}(\tau)F_{\mathrm{r}}(\tau)+X_{\mathrm{r}}(\tau)F_{\mathrm{i}}(\tau)-A_{\mathrm{i}}(\tau)X_{\mathrm{r}}(\tau)+A_{\mathrm{r}}(\tau)X_{\mathrm{i}}(\tau)\\=C_{\mathrm{i}}(\tau).	
	\end{matrix}\right.
	\end{split}
\end{equation}
Then, based on \eqref{eq.infer.realkroneckerproductvectorization}, \eqref{eq.dividetransccsme.variant} is finally converted to RFTVLSLE \eqref{eq.RFTVLSLE}.
\end{corollary}
Based on
\cite{He2024ZeroingND}, RFTVLSLE \eqref{eq.RFTVLSLE} is solved directly using the idea of Con-CZND2.
By means of the method of ZND
\cite{He2024ZeroingND}
in the real field, \eqref{eq.simplify.RFTVLSLE} is solved directly.

The error function is first defined as follows:
\begin{equation} \label{eq.define.errconcznd2}
	E_{\mathrm{R}}(\tau)=W_{\mathrm{R}}(\tau)X_{\mathrm{R}}(\tau)-B_{\mathrm{R}}(\tau),
\end{equation}
where $E_{\mathrm{R}}(\tau)\in\mathbb{R}^{2mn\times 1}$.
Next, the formula under the real field of ZND is employed to enable all elements of \eqref{eq.define.errconcznd2} 
to spontaneously converge to zero, which is obtained as
\begin{subequations} \label{eq.deduce.errconcznd2}
	\begin{align}
		\frac{\partial E_{\mathrm{R}}(\tau)}{\partial \tau} =-\gamma \psi \left ( E_{\mathrm{R}}(\tau) \right ),
		\\
		\dot{m}_{st}(\tau)=-\gamma \psi \left ( m_{st}(\tau) \right ).
	\end{align}
\end{subequations}
As in the previous subsection, $\gamma\in\mathbb{R^+}$ denotes the regulation parameter controlling the convergence rate, $\dot{m}_{st}(\tau)\in\mathbb{R}$ is $E_{\mathrm{R}}(\tau)$ elements differentiated from $\tau$,
and $\psi \left (\cdot  \right )$ denotes monotonically increasing odd activation function (MIOAF). For the subsequent study, a linear activation function $\psi \left (E_{\mathrm{R}}(\tau)  \right )=E_{\mathrm{R}}(\tau)$ is used in this case to eliminate the effect of MIOAF, and so \eqref{eq.deduce.errconcznd2} is simplified to:
\begin{subequations} \label{eq.infer.linearerrconcznd2}
	\begin{align}
		\frac{\partial E_{\mathrm{R}}(\tau)}{\partial \tau} =-\gamma E_{\mathrm{R}}(\tau),\\
		\dot{m}_{st}(\tau)=-\gamma m_{st}(\tau).
	\end{align}
\end{subequations}
Then, \eqref{eq.define.errconcznd2} is substituted into \eqref{eq.infer.linearerrconcznd2} to obtain \eqref{eq.join.linearerrconcznd2}:
\begin{equation} \label{eq.join.linearerrconcznd2}
	\begin{split}
		\dot{W}_{\mathrm{R}}(\tau)X_{\mathrm{R}}(\tau)+
		W_{\mathrm{R}}(\tau)\dot{X}_{\mathrm{R}}(\tau) -\dot{B}_{\mathrm{R}}(\tau)\\=-\gamma(W_{\mathrm{R}}(\tau)X_{\mathrm{R}}(\tau) -B_{\mathrm{R}}(\tau)).
	\end{split}
\end{equation}
Finally, the final solution model Con-CZND2 \eqref{eq.solve.linearerrconcznd2} is obtained:
\begin{equation} \label{eq.solve.linearerrconcznd2}
	\begin{split}
		\dot{X}_{\mathrm{R}}(\tau)
		=W^{+}_{\mathrm{R}}(\tau)(\dot{B}_{\mathrm{R}}(\tau)-\dot{W}_{\mathrm{R}}(\tau)X_{\mathrm{R}}(\tau)\\
		-\gamma(W_{\mathrm{R}}(\tau)X_{\mathrm{R}}(\tau) -B_{\mathrm{R}}(\tau))),
	\end{split}
\end{equation}
where $W^{+}_{\mathrm{R}}(\tau)$ is the pseudo-inverse matrix of $W_{\mathrm{R}}(\tau)$.

For ease of exposition, Example 3 in
\cite{He2024ZeroingND}
has been chosen. Symbols $s\left (\cdot  \right )$ and $c\left (\cdot  \right )$ denote the trigonometric functions $\sin\left (\cdot  \right )$ and $\cos\left (\cdot  \right )$, respectively. Then, 
\begin{equation} 
	F(\tau)
	=\begin{bmatrix}
		6+s(\tau)	&c(\tau) \\
		c(\tau)	&4+s(\tau)
	\end{bmatrix}
	+\mathrm{i}\begin{bmatrix}
		c(\tau)	& s(\tau) \\
		s(\tau) & c(\tau)
	\end{bmatrix}\in
	\mathbb{C}^{2\times 2},
	\notag
\end{equation}
\begin{equation}
	A(\tau)
	=\begin{bmatrix}
		c(\tau)	&s(\tau) \\
		-s(\tau)	&c(\tau)
	\end{bmatrix}
	+\mathrm{i}\begin{bmatrix}
		s(\tau)	& c(\tau) \\
		c(\tau) & -s(\tau)
	\end{bmatrix}\in
	\mathbb{C}^{2\times 2},
	\notag
\end{equation}
\begin{equation}
	C(\tau)
	=\begin{bmatrix}
		c_{\mathrm{r},11}(\tau) & c_{\mathrm{r},12}(\tau) \\
		c_{\mathrm{r},21}(\tau) & c_{\mathrm{r},22}(\tau) \\
	\end{bmatrix}
	+\mathrm{i}\begin{bmatrix}
		c_{\mathrm{i},11}(\tau) & 	c_{\mathrm{i},12}(\tau) \\
		c_{\mathrm{i},21}(\tau) & 	c_{\mathrm{i},22}(\tau) \\
	\end{bmatrix}\in
	\mathbb{C}^{2\times 2},
	\notag	
\end{equation}
where
$c_{\mathrm{r},11}(\tau)=2c^{2}(\tau)-2c(\tau)s(\tau)+6s(\tau)$,
$c_{\mathrm{r},12}(\tau)=4c(\tau)+2c(\tau)s(\tau)-2c^{2}(\tau)$, $c_{\mathrm{r},21}(\tau)=-2s(2\tau)-6c(\tau)+2$, 
$c_{\mathrm{r},22}(\tau)=2s(2\tau)-4s(\tau)-2$ 
and 
$c_{\mathrm{i},11}(\tau)=2c^{2}(\tau)+2c(\tau)s(\tau)+6s(\tau)$, $c_{\mathrm{i},12}(\tau)=4c(\tau)+2c(\tau)s(\tau)+2c^{2}(\tau)$, 
$c_{\mathrm{i},21}(\tau)=-2s(2\tau)-6c(\tau)-2$, $c_{\mathrm{i},22}(\tau)=-2s(2\tau)-4s(\tau)-2$.

\begin{figure*}[!t]\centering
	\subfloat[]{\includegraphics[width=1\columnwidth]{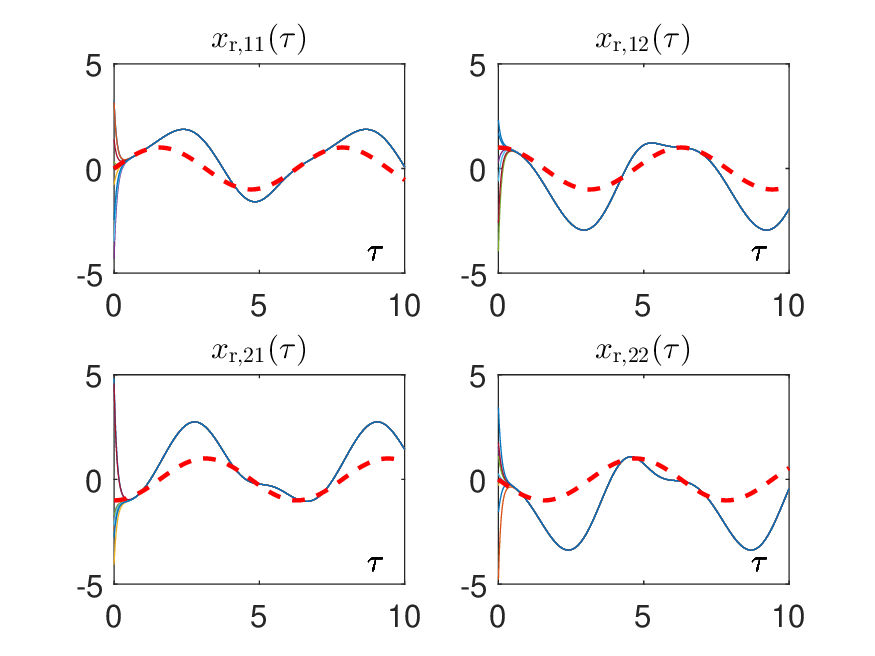}\label{fig.Con-CZND2.real}}
	\subfloat[]{\includegraphics[width=1\columnwidth]{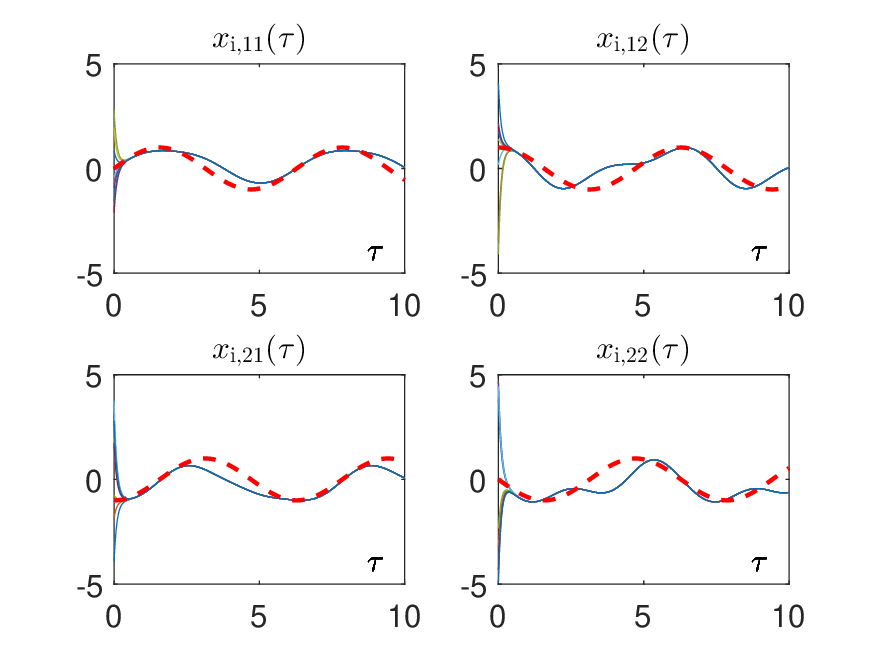}\label{fig.Con-CZND2.imaginary}}
	\caption{Solution $X_{\mathrm{R}}(\tau)$ computed by Con-CZND2 \eqref{eq.solve.linearerrconcznd2} model and exact solution $X_{\mathrm{R}}^*(\tau)$ ($\gamma$ equals 10).}
	\label{fig.Con-CZND2}
\end{figure*}

The unique exact solution to TVSSCME \eqref{eq.TVSSCME} is
\begin{equation}
	X^*(\tau)=\begin{bmatrix}
		s(\tau) & c(\tau)\\
		-c(\tau) & -s(\tau)\\
	\end{bmatrix}+\mathrm{i}\begin{bmatrix}
		s(\tau) & c(\tau) \\
		-c(\tau) & -s(\tau) \\
	\end{bmatrix}\in
	\mathbb{C}^{2\times 2},
\end{equation}
where $^{*}$ is the exact solution, same as below. It is deduced that the corresponding RETVLSLE \eqref{eq.RFTVLSLE} with the unique exact solution is:
\begin{equation}\label{eq.RETVLSLE.exact}
	X_{\mathrm{R}}^*(\tau)=\begin{bmatrix}
		\mathrm{vec}(X_{\mathrm{r}}^*(\tau))\\
		\mathrm{vec}(X_{\mathrm{i}}^*(\tau))	
	\end{bmatrix}=\begin{bmatrix}
		s(\tau) \\ -c(\tau)\\
		c(\tau) \\ -s(\tau)\\
		s(\tau) \\ -c(\tau) \\
		c(\tau) \\ -s(\tau) \\
	\end{bmatrix}\in
	\mathbb{R}^{8\times 1}.
\end{equation}
\begin{proof}
	According to 
	\cite{wuSolutionsMatrixEquations2006},
	for the time-invariant version of \eqref{eq.TVSSCME}, the condition for its unique solution is:
\begin{equation} \label{eq.SSCME.unique}
	\lambda (A\overline{A})\cap \lambda (F\overline{F})=\varnothing,
\end{equation}
where $\lambda$ is eigenvalue, $\cap$ is intersection, $\varnothing$ is empty set. For TVSSCME \eqref{eq.TVSSCME}, \eqref{eq.SSCME.unique} can also be generalized:
\begin{equation} \label{eq.TVSSCME.unique}
	\lambda (A(\tau)\overline{A}(\tau))\cap \lambda (F(\tau)\overline{F}(\tau))=\varnothing.
\end{equation}

When TVSSCME \eqref{eq.TVSSCME} is equivalently mapped to RFTVLSLE \eqref{eq.RFTVLSLE}, the judgment of unique \eqref{eq.RETVLSLE.exact} becomes the determinant $\mathrm{det}(\cdot)$ of the time-variant coefficients matrix $W_{\mathrm{R}}(\tau)$ is not equal to zero, i.e.,
\begin{equation} \label{eq.RFTVLSLE.unique}
	\mathrm{det}(W_{\mathrm{R}}(\tau))\ne 0.
\end{equation}

Then, the judgments of $\lambda (A(\tau)\overline{A}(\tau))\cap \lambda (F(\tau)\overline{F}(\tau))$ and $\mathrm{det}(W_{\mathrm{R}}(\tau))$ with time $\tau\in\left [ 0,10 \right ] $ are made simultaneously in Fig. \ref{fig.only.solve}. Figs. \ref{fig.only.solve}(a) and \ref{fig.only.solve}(b) are used to show whether the real and imaginary coefficients of $\lambda (A(\tau)\overline{A}(\tau))\cap \lambda (F(\tau)\overline{F}(\tau))$ will intersect at the same time $\tau$ point. Meanwhile, Fig. \ref{fig.only.solve}(c) with the same time $\tau$ is made to judge whether at some point $\mathrm{det}(W_{\mathrm{R}}(\tau))=0$.
\begin{figure*}[!t]\centering
	\subfloat[]{\includegraphics[width=0.6\columnwidth]{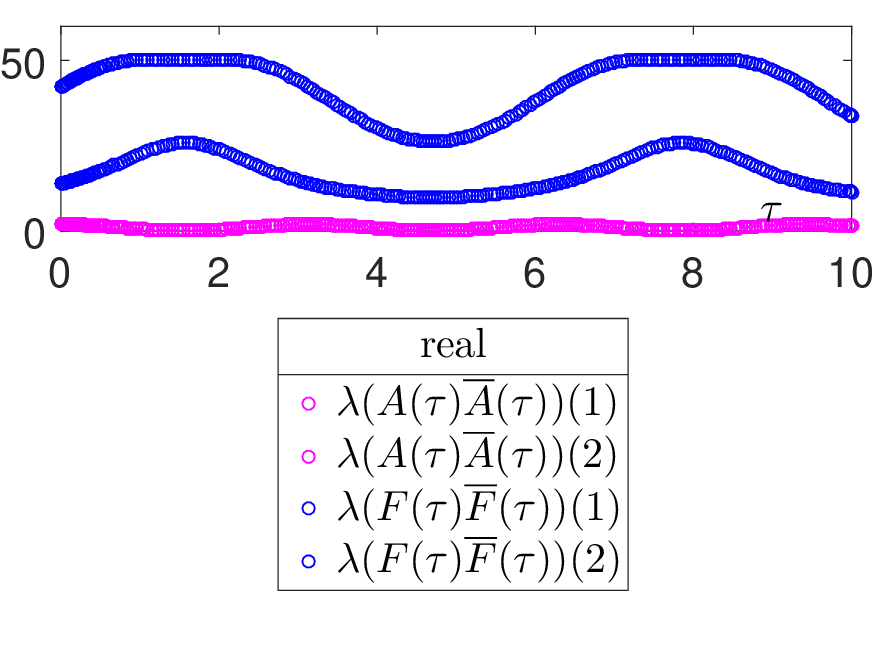}\label{fig.only.solve.real}}
	\subfloat[]{\includegraphics[width=0.6\columnwidth]{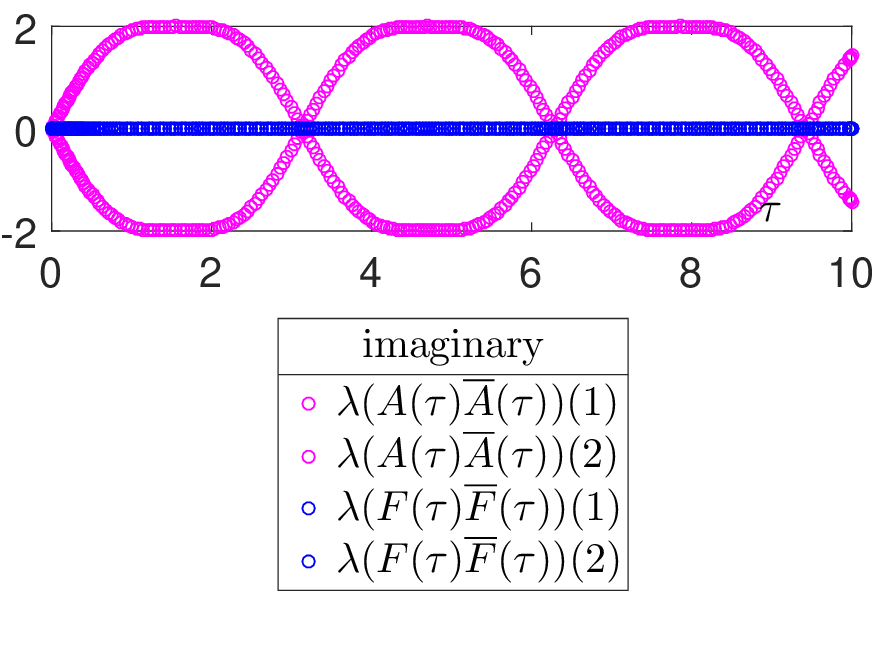}\label{fig.only.solve.imaginary}}
	\subfloat[]{\includegraphics[width=0.6\columnwidth]{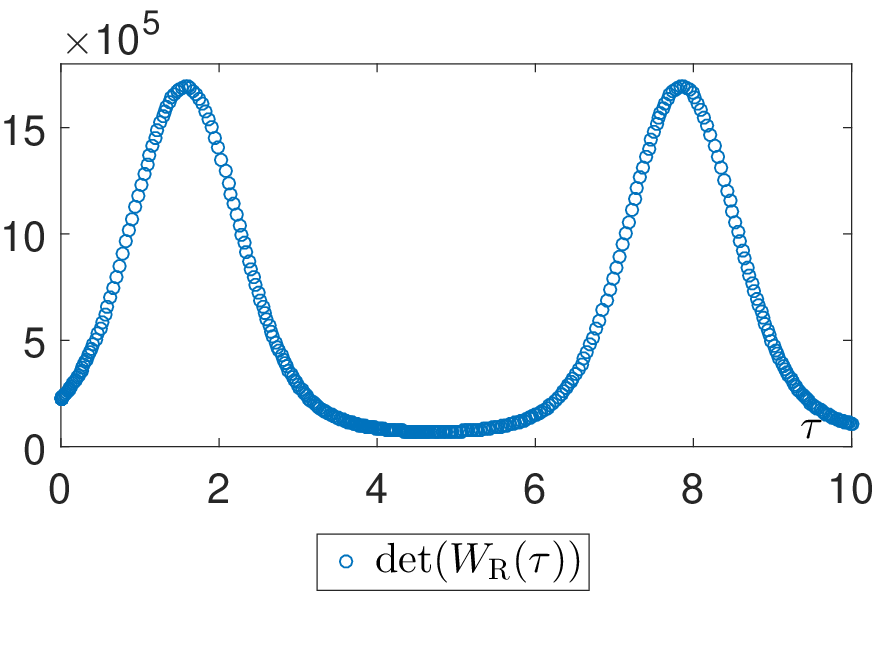}\label{fig.only.solve.det}}
	\caption{Judge whether $X^*(\tau)$ and $X_{\mathrm{R}}^*(\tau)$ are equivalent $\left (\lambda(\cdot)(n)\text{ denotes different eigenvalue for } n\in\mathbb{N^*}\right )$.}
	\label{fig.only.solve}
\end{figure*}

In Figs. \ref{fig.only.solve}(a) and \ref{fig.only.solve}(b), the real and imaginary coefficients between $\lambda (A(\tau)\overline{A}(\tau))$ and $\lambda (F(\tau)\overline{F}(\tau))$ do not intersect at the same time $\tau$ point. In Fig. \ref{fig.only.solve}(c), $\mathrm{det}(W_{\mathrm{R}}(\tau))\ne 0$. So $X^*(\tau)$ and $X_{\mathrm{R}}^*(\tau)$ are equivalent.

This proof is thus completed. 
\end{proof}
However, according to 
\cite{He2024ZeroingND},
since RFTVLSLE \eqref{eq.RFTVLSLE} is HD. Also, the main diagonal dominance of $W_{\mathrm{R}}(\tau)$ is not obvious. So the following experiment results are inconsistent with the theory that appears in Fig. \ref{fig.Con-CZND2}, where Con-CZND2 \eqref{eq.solve.linearerrconcznd2} model 
takes a random initial value $X_{0}\in\left [ -5,5 \right ] $, $\gamma$ equals 10, $\tau\in\left [ 0,10 \right ] $ by using ode45 function
\cite{He2024ZeroingND}.
The red dotted lines represent the exact solution for each solving element.

Thus, the solution by RFTVLSLE \eqref{eq.RFTVLSLE} error does not reach the exact solution. This is one of COD. However, in known studies, if conjugate operations occur, the solution using NN inevitably turns into RFTVLSLE \eqref{eq.RFTVLSLE}. So, again according to 
\cite{He2024ZeroingND}, it is found that the dimensional descent approach can be performed by an equivalent low-dimensional HN error matrix (which defaults to the complex matrix if not otherwise specified). And it is the source of Con-CZND1
\cite{He2024ZeroingND}
model.

Thus, a new dimensional descent approach is introduced next.
\section{Neural HN space compressive approximation approach}
According to 
\cite{He2024ZeroingND,dahlquistSpecialStabilityProblem1963},
reverse thinking is used here to reverse compress the system of RFTVLSLE \eqref{eq.RFTVLSLE} into equivalent TVSSCME \eqref{eq.TVSSCME} as shown in Fig. \ref{fig.NHNSCAA.matrix}.
\begin{figure}[!t]\centering
		\begin{tikzpicture}[  
		node distance=2cm, 
		auto, 
		every node/.style={align=center} 
		]  
		
		\node(node1) {  
			$\begin{aligned}  
				\begin{bmatrix}
					K_{11}(\tau)	&K_{12}(\tau) \\
					K_{21}(\tau)	&K_{22}(\tau)
				\end{bmatrix}
				\begin{bmatrix}
					\mathrm{vec}(X_{\mathrm{r}}(\tau))\\
					\mathrm{vec}(X_{\mathrm{i}}(\tau))	
				\end{bmatrix}
				=\begin{bmatrix}
					\mathrm{vec}(C_{\mathrm{r}}(\tau))\\
					\mathrm{vec}(C_{\mathrm{i}}(\tau))	
				\end{bmatrix} 
			\end{aligned}$  
		};  
		
		\node[below of=node1] (node2) {  
			$X(\tau)F(\tau)-A(\tau)\overline{X}(\tau)=C(\tau)$
		};  
		
		\draw[->] (node1) -- (node2);  
		
	\end{tikzpicture}
	 \caption{Schematic of neural HN space compressive approximation approach (NHNSCAA) operation by reverse thinking.}
	 \label{fig.NHNSCAA.matrix}
\end{figure}

This is proposed next: neural HN space compressive approximation approach (NHNSCAA).
\begin{definition}
Neural HN space compressive approximation approach (NHNSCAA) is defined as follows: 

If a time-variant large-scale linear equation error occurs with high dimension, this error can be replaced by mapping it to an equivalent low-dimensional time-variant hypercomplex number matrix equation error. The analytical solution can then be further approximated to the theoretical solution.
\end{definition}
For the sake of rigor, the default here is that both time-variant LSLE and the low-dimensional HN matrix equations have only equivalent unique solutions, like RFTVLSLE \eqref{eq.RFTVLSLE} and TVSSCME \eqref{eq.TVSSCME}. At the same time,
\begin{equation} \label{eq.define.errconcznd1}
	E_{\mathrm{C1}}(\tau)=X(\tau)F(\tau)-A(\tau)\overline{X}(\tau)-C(\tau),
\end{equation}
where $E_{\mathrm{C1}}(\tau)\in\mathbb{C}^{m\times n}$ is the compressed HN matrix error. The number of elements in $E_{\mathrm{C1}}(\tau)$ is only $m\times n$. However, in Fig. \ref{fig.complexmatrix}, $E_{\mathrm{C1}}(\tau)$ is essentially a $ \mathbb{R}^{2mn\times 1}\mapsto \mathbb{C}^{m\times n}$ compression. $\mathbb{C}^{m\times n}$ can be regarded as $\mathbb{R}^{m\times n \times 2}$. $ \mathbb{R}^{m\times n \times 2} $ represents a tensor with 3 measures, which is relevant to deep learning
\cite{LIAO2022440}.
Therefore, $ \mathbb{R}^{2mn\times 1}\mapsto \mathbb{C}^{m\times n}$ can also be thought of as tensor optimization
\cite{Kolda2009TensorDA}.

\begin{figure*}[!t]  
	\centering  
	\subfloat[]{%
		
		\resizebox{0.45\textwidth}{!}{%
			\def\layersep{2cm}
			\raisebox{3cm}{
				\begin{tikzpicture}[shorten >=1pt,->,draw=black!50, node distance=\layersep]
					\tikzstyle{every pin edge}=[<-,shorten <=1pt]
					\tikzstyle{neuron}=[circle,fill=black!25,minimum size=17pt,inner sep=0pt]
					\tikzstyle{input neuron}=[neuron, fill=green!50];
					\tikzstyle{output neuron}=[neuron, fill=red!50];
					\tikzstyle{hidden neuron12}=[neuron, fill=white!50];
					\tikzstyle{hidden neuron2}=[neuron, fill=green!50];
					\tikzstyle{hidden neuron4}=[neuron, fill=red!50];
					\tikzstyle{annot} = [text width=4em, text centered]
					
					\foreach \name / \y in {1,...,8}
					\node[input neuron, pin=left:Random Input \y] (I-\name) at (0,-\y) {};

					\foreach \name / \y in {1}
					\path[yshift=-9.5cm]
					node[hidden neuron12,right of=I] (H12-\name) at (0.5*\layersep,-\y cm){...};
					
					\foreach \name / \y in {1,...,8}
					\path[yshift=0cm]
					node[hidden neuron2,right of=H2] (H2-\name) at (2*\layersep,-\y cm){};
					
					\foreach \name / \y in {1,...,8}
					\path[yshift=0cm]
					node[hidden neuron4,right of=H3,pin={[pin edge={->}]right:Output},] (H3-\name) at (3*\layersep,-\y cm){};
					
					\foreach \source in {1,...,8}
					\path (I-\source) edge (H12-1);

					\foreach \dest in {1,...,8}
					\path (H12-1) edge (H2-\dest);

					\foreach \source in {1,...,8}
					\foreach \dest in {1,...,8}
					\path (H2-\source) edge (I-\dest);
					
					\foreach \source in {1,...,8}		
					\path (H2-\source) edge (H3-\source);
					
					\node[annot,above of=I-1, node distance=1cm] (hl) {Input Real Layer};
					\node[annot,right of=hl] (h2){};
					\node[annot,right of=h2] (h3){};
					\node[annot,right of=h3] (h4){Hidden Real Layer1};
					\node[annot,right of=h4] {Output Real Layer};
			\end{tikzpicture}}
		}%
		
		\label{fig.NHNSCAA.before}
	} 
	\hfill 
	\subfloat[]{%
		\resizebox{0.45\textwidth}{!}{%
			\def\layersep{2cm}
			
			\begin{tikzpicture}[shorten >=1pt,->,draw=black!50, node distance=\layersep]
				\tikzstyle{every pin edge}=[<-,shorten <=1pt]
				\tikzstyle{neuron}=[circle,fill=black!25,minimum size=17pt,inner sep=0pt]
				\tikzstyle{input neuron}=[neuron, fill=green!50];
				\tikzstyle{output neuron}=[neuron, fill=red!50];
				
				\tikzstyle{hidden neuron1}=[neuron, fill=blue!50];
				\tikzstyle{hidden neuron2}=[neuron, fill=blue!50];
				\tikzstyle{hidden neuron23}=[neuron, fill=white!50];
				\tikzstyle{hidden neuron3}=[neuron, fill=green!50];
				\tikzstyle{hidden neuron4}=[neuron, fill=red!50];
				\tikzstyle{annot} = [text width=4em, text centered]
				
				\foreach \name / \y in {1,...,8}
				\node[input neuron, pin=left:Random Input \y] (I-\name) at (0,-\y) {};
				
				\foreach \name / \y in {1,...,4}
				\path[yshift=-8cm]
				node[hidden neuron1] (H1-\name) at (\layersep,-\y cm) {};
				\foreach \name / \y in {1,...,4}
				\path[yshift=-8cm]
				node[hidden neuron2,right of=H1] (H2-\name) at (\layersep,-\y cm){};
				
				\foreach \name / \y in {1}
				\path[yshift=-9.5cm]
				node[hidden neuron23,right of=H2] (H23-\name) at (0.5*\layersep,-\y cm){...};
				
				\foreach \name / \y in {1,...,8}
				\path[yshift=0cm]
				node[hidden neuron3,right of=H2] (H3-\name) at (2*\layersep,-\y cm){};
				
				\foreach \name / \y in {1,...,8}
				\path[yshift=0cm]
				node[hidden neuron4,right of=H3,pin={[pin edge={->}]right:Output},] (H4-\name) at (3*\layersep,-\y cm){};
				
				\foreach \source in {1,...,4}
				\path (I-\source) edge (H1-\source);
				\path (I-5) edge (H1-1);
				\path (I-6) edge (H1-2);
				\path (I-7) edge (H1-3);
				\path (I-8) edge (H1-4);
				
				\foreach \source in {1,...,4}
				\path (H1-\source) edge (H23-1);
				
				\foreach \dest in {1,...,4}
				\path (H23-1) edge (H2-\dest);
				
				\foreach \source in {1,...,4}
				\path (H2-\source) edge (H3-\source);
				\path (H2-1) edge (H3-5);
				\path (H2-2) edge (H3-6);
				\path (H2-3) edge (H3-7);
				\path (H2-4) edge (H3-8);
				
				\foreach \source in {1,...,8}
				\path (H3-\source) edge (H4-\source);
				
				\foreach \source in {1,...,8}
				\foreach \dest in {1,...,8}
				\path (H3-\source) edge (I-\dest);
				
				
				
				\node[annot,below of=H1-1, node distance=4cm] (hl) {Hidden Complex Layer1};
				\node[annot,above of=I-1, node distance=1cm](IRltext) {Input Real Layer};
				\node[annot,right of=hl] (h2){Hidden Complex Layer2};
				\node[annot,right of=IRltext, node distance=6cm] (h3){Hidden Real Layer1};
				\node[annot,right of=h3] {Output Real Layer};
			\end{tikzpicture}
		}%
		\label{fig.NHNSCAA.after}
	}  
	\caption{Schematic of neural HN space compressive approximation approach (NHNSCAA) mechanism. (a) Before processing NHNSCAA. (b) After processing NHNSCAA.} 
	\label{fig.NHNSCAA.network} 
\end{figure*}
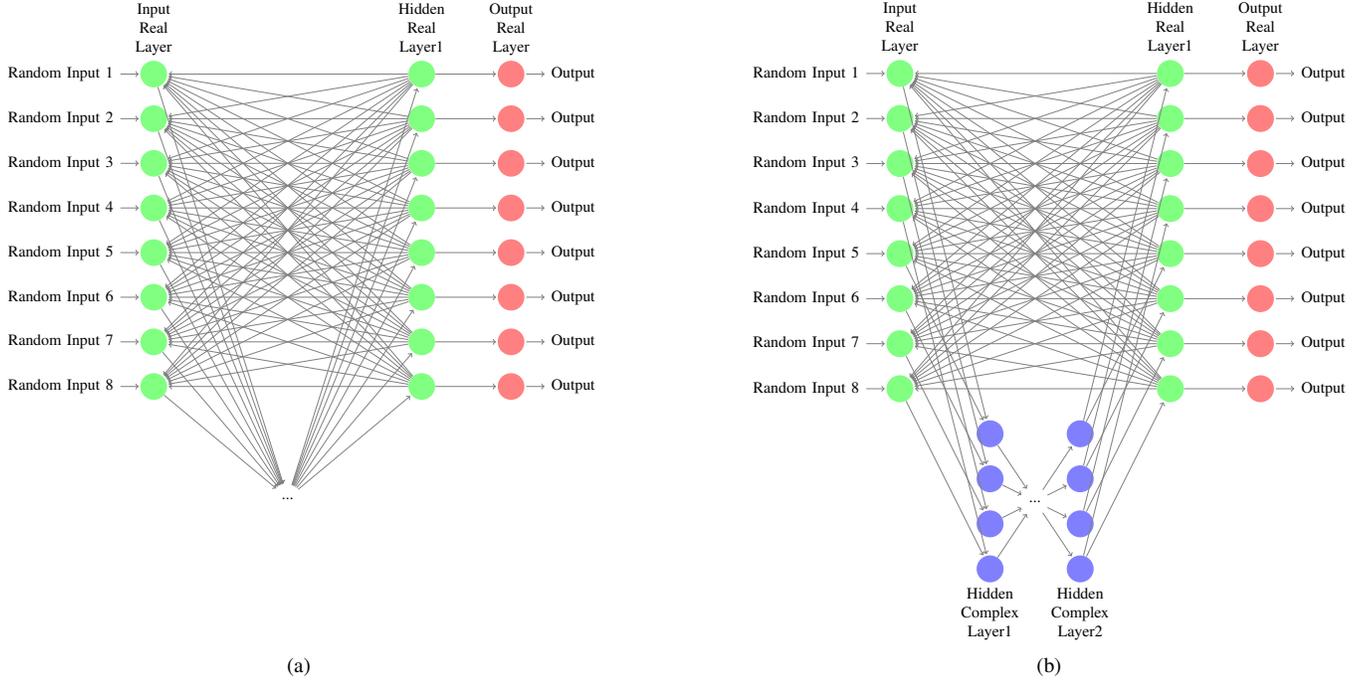

However, for the sake of mathematical rigor, the tensor of multiple measures is uniformly referred to as HN matrix if not otherwise specified. The structure of the optimized NN is shown in Fig. \ref{fig.NHNSCAA.network}.

According to the idea of NHNSCAA, \eqref{eq.TVSSCME.conj} is obtained by taking the conjugate for both sides of \eqref{eq.TVSSCME} simultaneously.
\begin{equation} \label{eq.TVSSCME.conj}
	\begin{split}
		\overline{X}(\tau)\overline{F}(\tau)-\overline{A}(\tau)X(\tau)=\overline{C}(\tau),
	\end{split}
\end{equation}
where \eqref{eq.TVSSCME.conj} also have equivalent unique solution with \eqref{eq.RFTVLSLE}.
\begin{proof}
Like \eqref{eq.transccsme.variant}, \eqref{eq.TVSSCME.conj} can be split into \eqref{eq.transccsme.variant.conj} as follows.
\begin{equation}\label{eq.transccsme.variant.conj}
	\begin{split}
		(X_{\mathrm{r}}(\tau)-\mathrm{i}X_{\mathrm{i}}(\tau))(F_{\mathrm{r}}(\tau)-\mathrm{i}F_{\mathrm{i}}(\tau))\\-(A_{\mathrm{r}}(\tau)-\mathrm{i}A_{\mathrm{i}}(\tau))(X_{\mathrm{r}}(\tau)+\mathrm{i}X_{\mathrm{i}}(\tau))\\=(C_{\mathrm{r}}(\tau)-\mathrm{i}C_{\mathrm{i}}(\tau)).
	\end{split}	
\end{equation}

Separating \eqref{eq.transccsme.variant.conj} is then carried out to obtain the equivalent real field \eqref{eq.dividetransccsme.variant.conj} as follows.
\begin{equation}\label{eq.dividetransccsme.variant.conj}
	\begin{split}
		\left\{\begin{matrix}
			X_{\mathrm{r}}(\tau)F_{\mathrm{r}}(\tau)-X_{\mathrm{i}}(\tau)F_{\mathrm{i}}(\tau)-A_{\mathrm{r}}(\tau)X_{\mathrm{r}}(\tau)-A_{\mathrm{i}}(\tau)X_{\mathrm{i}}(\tau)\\=C_{\mathrm{r}}(\tau),	\\
			-( X_{\mathrm{i}}(\tau)F_{\mathrm{r}}(\tau)+X_{\mathrm{r}}(\tau)F_{\mathrm{i}}(\tau)-A_{\mathrm{i}}(\tau)X_{\mathrm{r}}(\tau)+A_{\mathrm{r}}(\tau)X_{\mathrm{i}}(\tau))\\=-C_{\mathrm{i}}(\tau).	
		\end{matrix}\right.
	\end{split}
\end{equation}
Notice that \eqref{eq.dividetransccsme.variant.conj} is equivalent to \eqref{eq.dividetransccsme.variant}. And \eqref{eq.TVSSCME.unique} is applied to \eqref{eq.TVSSCME.conj} as well.

This proof is thus completed. 	
\end{proof}
The following error matrix of \eqref{eq.TVSSCME.conj} is then defined:
\begin{equation} \label{eq.define.errconcznd1.conj}
	E_{\mathrm{C2}}(\tau)=\overline{X}(\tau)\overline{F}(\tau)-\overline{A}(\tau)X(\tau)-\overline{C}(\tau),
\end{equation}
where $E_{\mathrm{C2}}(\tau)\in\mathbb{C}^{m\times n}$ is the compressed HN matrix error. The number of elements of $E_{\mathrm{C2}}(\tau)$ is only $m\times n$.

Next, a method similar to Con-CZND1 
\cite{He2024ZeroingND}
is employed to enable all elements of \eqref{eq.define.errconcznd1.conj} to spontaneously converge to zero, i.e.,
\begin{subequations} \label{eq.deduce.errconcznd1.conj}
	\begin{align}
		\frac{\partial E_{\mathrm{C2}}(\tau)}{\partial \tau} =-\gamma \psi \left ( E_{\mathrm{C2}}(\tau) \right ),
		\\
		\dot{\tilde{m}}_{st}(\tau)=-\gamma \psi \left ( \tilde{m}_{st}(\tau) \right ).
	\end{align}
\end{subequations}
As in the previous subsection, $\gamma\in\mathbb{R^+}$ denotes the regulation parameter controlling the convergence rate, $\dot{\tilde{m}}_{st}(\tau)\in\mathbb{C}$ is $E_{\mathrm{C2}}(\tau)$ elements differentiated from $\tau$,
and $\psi \left (\cdot  \right )$ denotes MIOAF. For the subsequent study, a linear activation function $\psi \left (E_{\mathrm{C2}}(\tau)  \right )=E_{\mathrm{C2}}(\tau)$ is used in this case to eliminate the effect of MIOAF, and so \eqref{eq.deduce.errconcznd1.conj} is simplified to:
\begin{subequations} \label{eq.infer.linearerrconcznd1.conj}
	\begin{align}
		\frac{\partial E_{\mathrm{C2}}(\tau)}{\partial \tau} =-\gamma E_{\mathrm{C2}}(\tau),\\
		\dot{\tilde{m}}_{st}(\tau)=-\gamma  \tilde{m}_{st}(\tau).
	\end{align}
\end{subequations}
Then, \eqref{eq.define.errconcznd1.conj} is substituted into \eqref{eq.infer.linearerrconcznd1.conj} to obtain \eqref{eq.join.linearerrconcznd1.conj}:
\begin{equation} \label{eq.join.linearerrconcznd1.conj}
	\begin{split}
	\dot{\overline{X}}(\tau)\overline{F}(\tau)+\overline{X}(\tau)\dot{\overline{F}}(\tau)-\dot{\overline{A}}(\tau)X(\tau)-\overline{A}(\tau)\dot{X}(\tau)-\dot{\overline{C}}(\tau)
	\\
	=-\gamma(\overline{X}(\tau)\overline{F}(\tau)-\overline{A}(\tau)X(\tau)-\overline{C}(\tau)).\end{split}
\end{equation}
Using \eqref{eq.infer.complexkroneckerproductvectorization}, \eqref{eq.join.linearerrconcznd1.conj} is converted to \eqref{eq.use.complexkroneckerproductvectorization.conj}:
\begin{equation} \label{eq.use.complexkroneckerproductvectorization.conj}
	\begin{split}
		(F^{\mathrm{H}}(\tau)\otimes I_{m})\mathrm{vec}(\dot{\overline{X}}(\tau))-(\overline{I_{n}^{\mathrm{H}}}\otimes \overline{A}(\tau))\mathrm{vec}(\dot{X}(\tau))
		\\
		=\mathrm{vec}(\dot{\overline{C}}(\tau)+\dot{\overline{A}}(\tau)X(\tau)-\overline{X}(\tau)\dot{\overline{F}}(\tau))
		\\
		-\gamma\mathrm{vec}(\overline{X}(\tau)\overline{F}(\tau)-\overline{A}(\tau)X(\tau)-\overline{C}(\tau)).
	\end{split}	
\end{equation}
Then, \eqref{eq.use.complexkroneckerproductvectorization.conj} is further reformulated as
\begin{equation}
	\label{eq.simplify.complexkroneckerproductvectorization.conj}
	H(\tau)\mathrm{vec}(\dot{\overline{X}}(\tau))-L(\tau)\mathrm{vec}(\dot{X}(\tau))
	=O(\tau),
\end{equation}
where $H(\tau)=(F^{\mathrm{H}}(\tau) \otimes I_{m})\in
\mathbb{C}^{nm\times mn}$, $L(\tau)=(\overline{I_{n}^{\mathrm{H}}} \otimes \overline{A}(\tau))=(I_{n}\otimes \overline{A}(\tau))\in
\mathbb{C}^{mn\times nm}$,
$O(\tau)=\mathrm{vec}(\dot{\overline{C}}(\tau)+\dot{\overline{A}}(\tau)X(\tau)-\overline{X}(\tau)\dot{\overline{F}}(\tau))
-\gamma\mathrm{vec}(\overline{X}(\tau)\overline{F}(\tau)-\overline{A}(\tau)X(\tau)-\overline{C}(\tau))\in
\mathbb{C}^{mn\times 1}$. Based on the linearity of the derivative as well as \eqref{eq.define.complexmatrix}, \eqref{eq.simplify.complexkroneckerproductvectorization.conj} can be written in the form of the following real-only matrix operation:
\begin{equation}\label{eq.divide.complexkroneckerproductvectorization.conj}
	\begin{bmatrix}
		H_{\mathrm{r}}(\tau)-L_{\mathrm{r}}(\tau)	&H_{\mathrm{i}}(\tau)+L_{\mathrm{i}}(\tau) \\
		H_{\mathrm{i}}(\tau)-L_{\mathrm{i}}(\tau)	&-(H_{\mathrm{r}}(\tau)+L_{\mathrm{r}}(\tau))
	\end{bmatrix}
	\begin{bmatrix}
		\dot{Z}_{\mathrm{r}}(\tau)\\
		\dot{Z}_{\mathrm{i}}(\tau)	
	\end{bmatrix}
	=\begin{bmatrix}
		O_{\mathrm{r}}(\tau)\\
		O_{\mathrm{i}}(\tau)	
	\end{bmatrix},
\end{equation}
where $Z(\tau)=\mathrm{vec}(X(\tau))\in
\mathbb{C}^{mn\times 1}$, $\dot{Z}(\tau)=\mathrm{vec}(\dot{X}(\tau))\in
\\\mathbb{C}^{mn\times 1}$. To simplify, let $W_{\mathrm{C2}}(\tau)=\left[H_{\mathrm{r}}(\tau)-L_{\mathrm{r}}(\tau), H_{\mathrm{i}}(\tau)+\right.\\\left.L_{\mathrm{i}}(\tau);
H_{\mathrm{i}}(\tau)-L_{\mathrm{i}}(\tau),-(H_{\mathrm{r}}(\tau)+L_{\mathrm{r}}(\tau))\right ]
\in\mathbb{R}^{2mn\times 2mn}$, $\dot{X}_{\mathrm{C2}}(\tau)= \left [\dot{Z}_{\mathrm{r}}(\tau);
\dot{Z}_{\mathrm{i}}(\tau)\right ]\in
\mathbb{R}^{2mn\times 1}$, $B_{\mathrm{C2}}(\tau)= \left [O_{\mathrm{r}}(\tau);
O_{\mathrm{i}}\right.\\\left.(\tau)\right ]\in\mathbb{R}^{2mn\times 1}$. The final solution model Con-CZND1\_conj is obtained:
\begin{equation} \label{eq.solve.linearerrconcznd1.conj}
	\dot{X}_{\mathrm{C2}}(\tau) =W^{+}_{\mathrm{C2}}(\tau)B_{\mathrm{C2}}(\tau),
\end{equation}
where $W^{+}_{\mathrm{C2}}(\tau)$ is the pseudo-inverse matrix of $W_{\mathrm{C2}}(\tau)$.
\begin{theorem}\label{thm.1}
Given differentiable time-variant matrices $F(\tau)\in
\mathbb{C}^{n\times n}$, $A(\tau)\in
\mathbb{C}^{m\times m}$, and $C(\tau)\in
\mathbb{C}^{m\times n}$, if RETVLSLE \eqref{eq.RFTVLSLE} only has one theoretical time-variant solution $X_{\mathrm{R}}^*(\tau)\in
\mathbb{R}^{2mn\times 1}$, then each solving element of \eqref{eq.solve.linearerrconcznd1.conj} converges to the corresponding theoretical time-variant solving elements.
	\begin{proof}
	Omitted because the proof process is similar to Theorem 2 of \cite{He2024ZeroingND}.	
	\end{proof}
\end{theorem}
\section{Numerical Experimentation, Verification and Extension}
The templates are intended to highlight the validity of Con-CZND1\_conj model under numerical Example 3 using the same 
\cite{He2024ZeroingND}.
Also, $\gamma$ - swallowed phenomenon is introduced which is not mentioned in known studies. This is the biggest difference toward space compressive approximation to activate function
and sampling discretion.
\begin{remark}
	Unless otherwise stated, the following numerical experiments all take a random initial value $X_{0}\in\left [ -5,5 \right ] $, $\tau\in\left [ 0,10 \right ] $ by using ode45
	\cite{He2024ZeroingND}
	function. The red dotted lines represent the exact solution for each solving element. Because of RFTVLSLE \eqref{eq.RFTVLSLE}, the residuals are uniformly defined $\left \|X_{\mathrm{R}}(\tau)-X_{\mathrm{R}}^*(\tau)   \right \|_{\mathrm{F}}$, where $\left \|\cdot   \right \|_{\mathrm{F}}$ stands for Frobenius norm.
\end{remark}
\subsection{Con-CZND1\_conj model}
\begin{figure*}[!t]\centering
	\subfloat[]{\includegraphics[width=1\columnwidth]{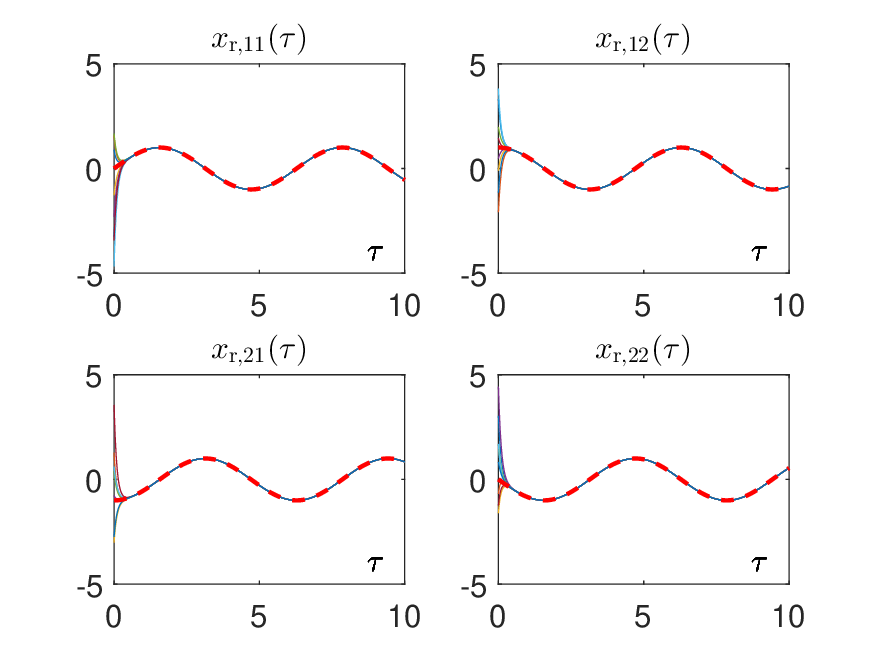}\label{fig.Con-CZND1_conj.real}}
	\subfloat[]{\includegraphics[width=1\columnwidth]{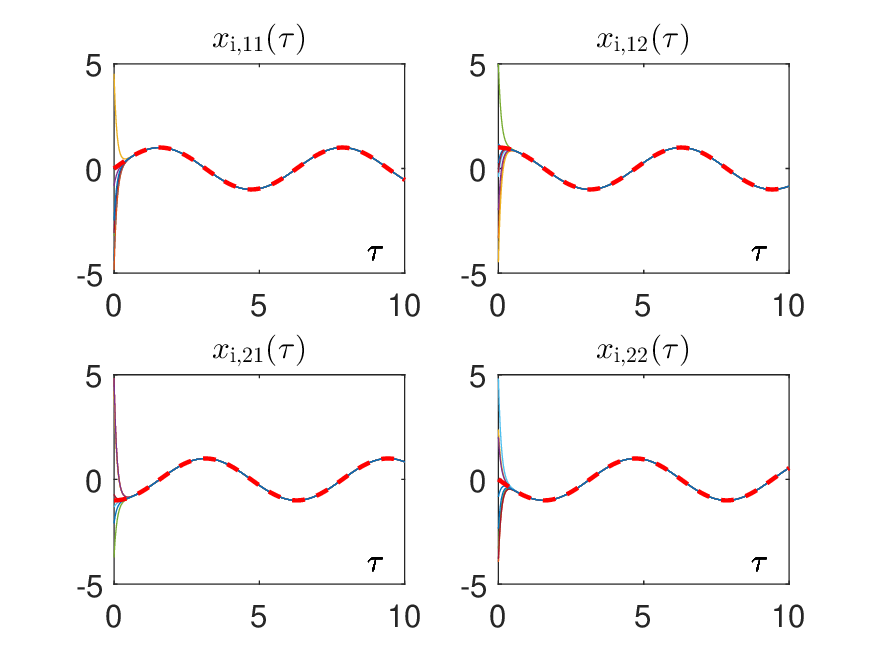}\label{fig.Con-CZND1_conj.imaginary}}
	\caption{Solution $X_{\mathrm{R}}(\tau)$ computed by Con-CZND1\_conj \eqref{eq.solve.linearerrconcznd1.conj} model and exact solution $X_{\mathrm{R}}^*(\tau)$ ($\gamma$ equals 10).}
	\label{fig.Con-CZND1_conj}
\end{figure*}

\begin{figure*}[!t]\centering
	\subfloat[]{\includegraphics[width=1\columnwidth]{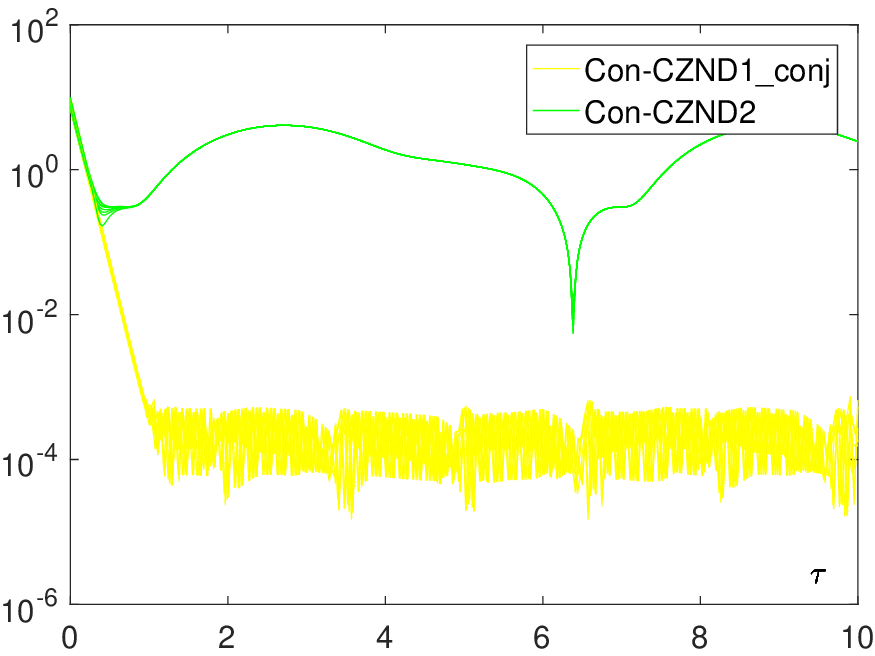}\label{fig.Con-CZND1_conj.vs.Con-CZND2}}
	\subfloat[]{\includegraphics[width=1\columnwidth]{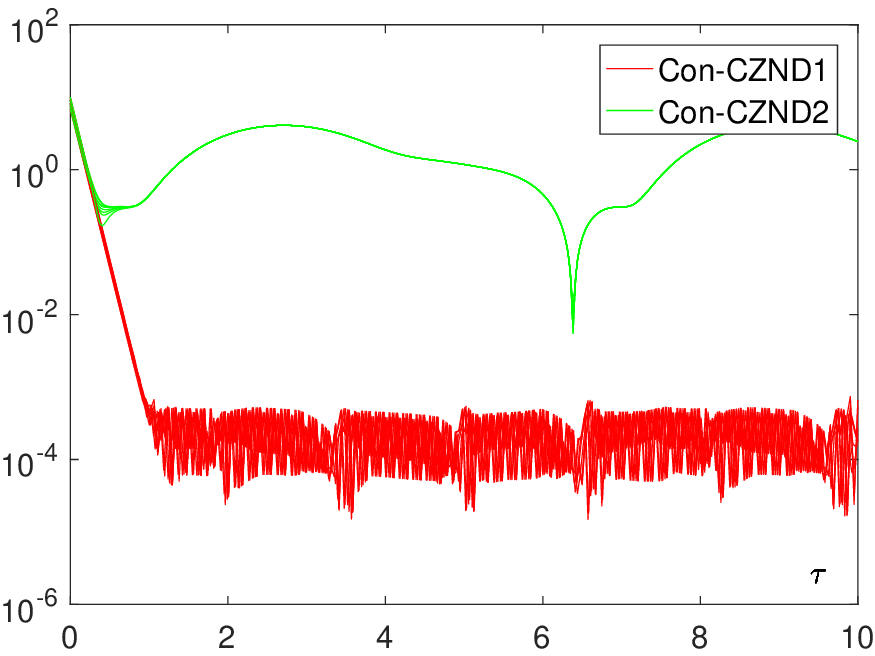}\label{fig.Con-CZND1.vs.Con-CZND2}}
	\caption{Logarithmic residual $\left \|X_{\mathrm{R}}(\tau)-X_{\mathrm{R}}^*(\tau)   \right \|_{\mathrm{F}}$ trajectories computed ($\gamma$ equals 10). (a) Con-CZND1\_conj \eqref{eq.solve.linearerrconcznd1.conj} model vs. Con-CZND2 \eqref{eq.solve.linearerrconcznd2} model. (b) Con-CZND1\cite{He2024ZeroingND} model vs. Con-CZND2 \eqref{eq.solve.linearerrconcznd2} model.}
	\label{fig.Con-CZND1_conj and fig.Con-CZND1 vs fig.Con-CZND2}
\end{figure*}

Executing Con-CZND1\_conj \eqref{eq.solve.linearerrconcznd1.conj} model when $\gamma$ equals 10, the results of the numerical experiments are shown in Fig. \ref{fig.Con-CZND1_conj}. The logarithmic residual $\left \|X_{\mathrm{R}}(\tau)-X_{\mathrm{R}}^*(\tau)   \right \|_{\mathrm{F}}$ trajectories of Con-CZND1\_conj \eqref{eq.solve.linearerrconcznd1.conj} model vs Con-CZND2 \eqref{eq.solve.linearerrconcznd2} model is given in Fig. \ref{fig.Con-CZND1_conj and fig.Con-CZND1 vs fig.Con-CZND2}(a), with attaching Con-CZND1
\cite{He2024ZeroingND}
model vs Con-CZND2 \eqref{eq.solve.linearerrconcznd1.conj} model in Fig. \ref{fig.Con-CZND1_conj and fig.Con-CZND1 vs fig.Con-CZND2}(b).

As can be seen from Fig. \ref{fig.Con-CZND1_conj and fig.Con-CZND1 vs fig.Con-CZND2}, Con-CZND1\_conj \eqref{eq.solve.linearerrconcznd1.conj} model
and Con-CZND1 model have consistent residual curves since $E_{\mathrm{C2}}(\tau)\in\mathbb{C}^{m\times n}$ and $E_{\mathrm{C1}}(\tau)\in\mathbb{C}^{m\times n}$. Therefore, it can be consistently determined that HN matrix equation \eqref{eq.TVSSCME.conj}, which is compressed in terms of RFTVLSLE \eqref{eq.RFTVLSLE}, only exchanges space with TVSSCME \eqref{eq.TVSSCME} by the real $\gamma$. There is no significant difference between HN matrix equation \eqref{eq.TVSSCME.conj} and TVSSCME \eqref{eq.TVSSCME}. Next, $\gamma$ will do some HN extension and the corresponding numerical experiments will be supplied to highlight NHNSCAA.
\subsection{$\gamma$ - swallowed and Space compressive approximation}
\begin{figure*}[!t]\centering
	\subfloat[]{\includegraphics[width=1\columnwidth]{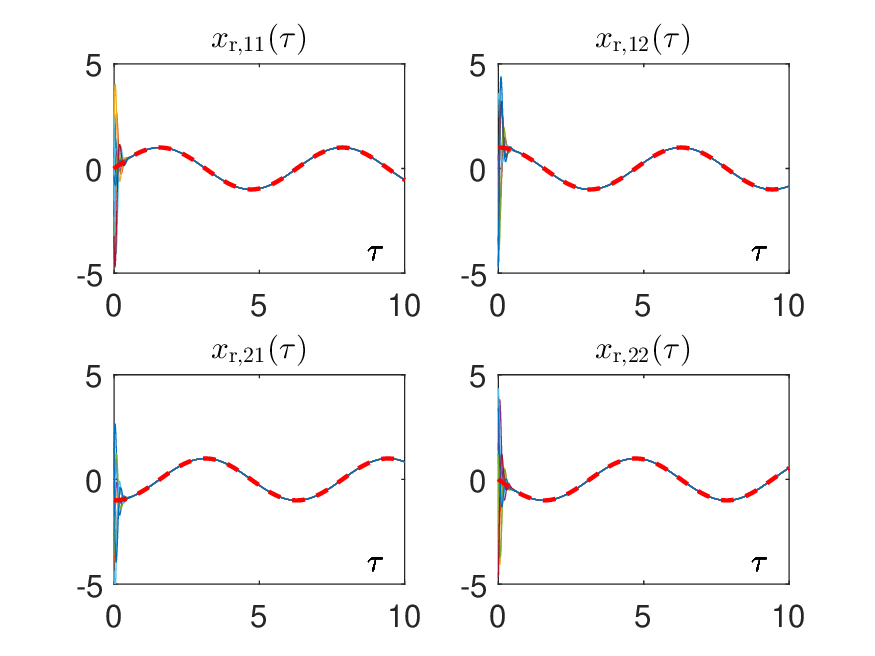}\label{fig.Con-CZND1_conj.real.extend10+20i}}
	\subfloat[]{\includegraphics[width=1\columnwidth]{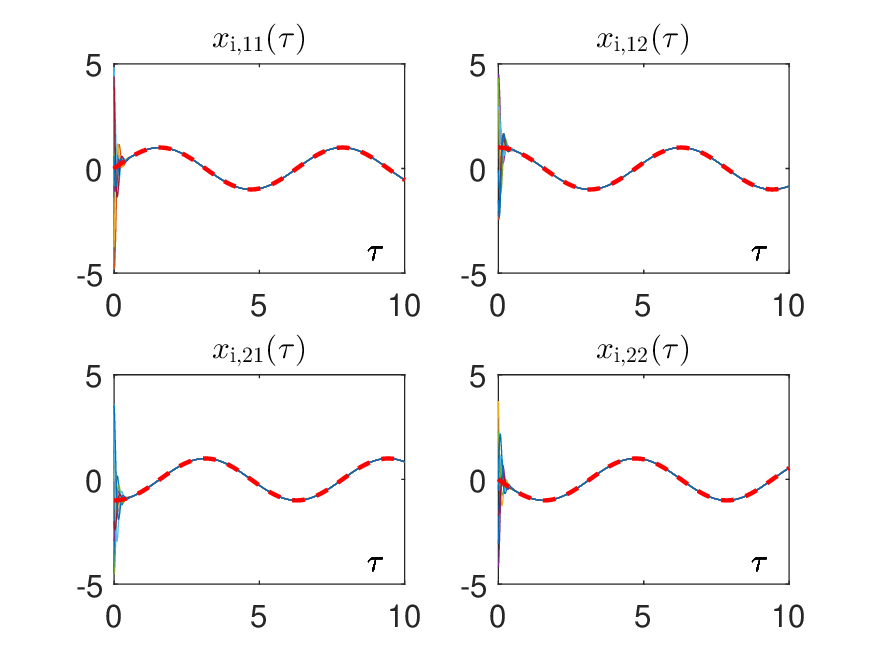}\label{fig.Con-CZND1_conj.imaginary.extend10+20i}}
	\caption{Solution $X_{\mathrm{R}}(\tau)$ computed by Con-CZND1\_conj \eqref{eq.solve.linearerrconcznd1.conj} model and exact solution $X_{\mathrm{R}}^*(\tau)$ ($\gamma$ equals 10$+$20$\mathrm{i}$).}
	\label{fig.Con-CZND1_conj.extend10+20i}
\end{figure*}

\begin{figure*}[!t]\centering
	\subfloat[]{\includegraphics[width=1\columnwidth]{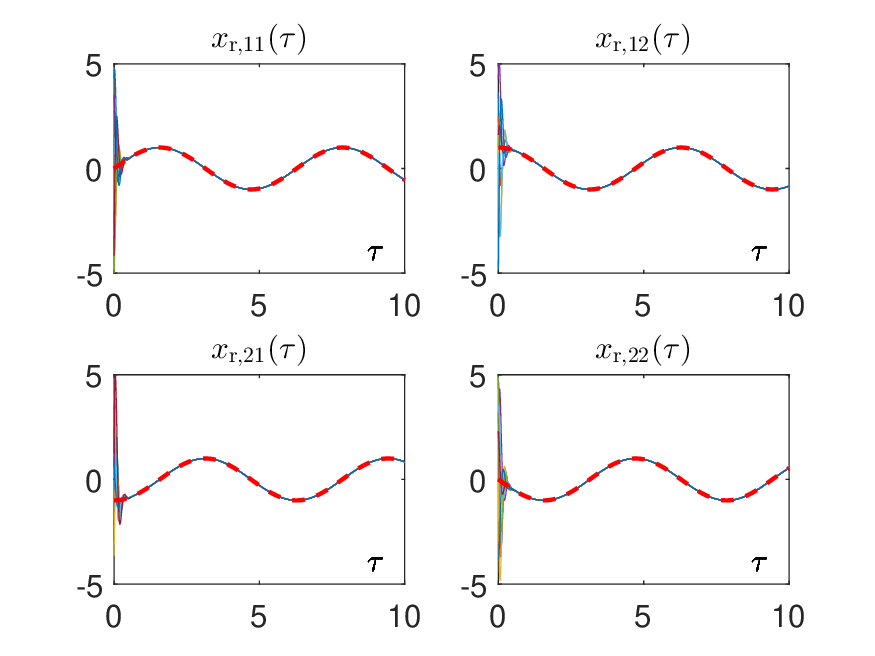}\label{fig.Con-CZND1.real.extend10+20i}}
	\subfloat[]{\includegraphics[width=1\columnwidth]{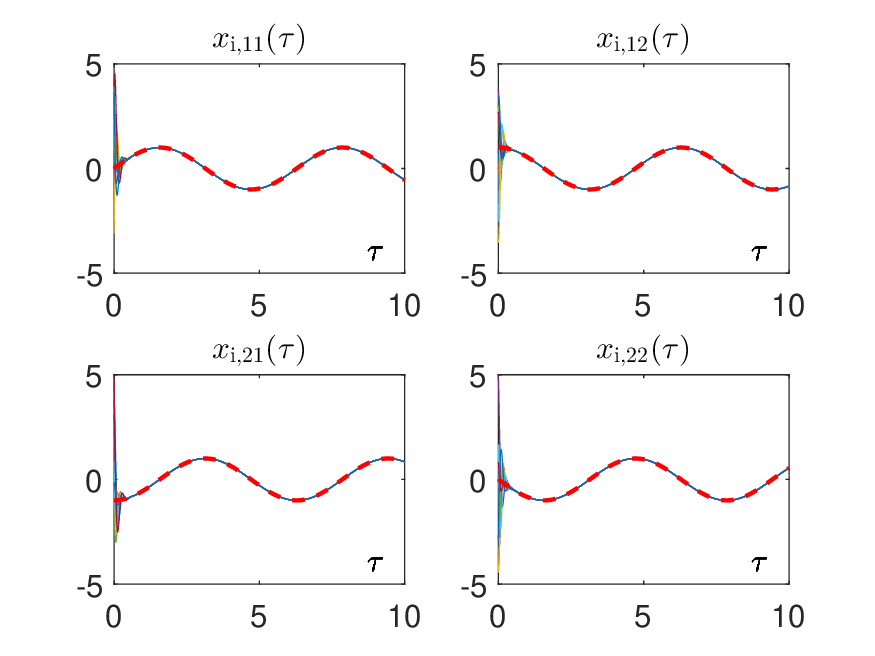}\label{fig.Con-CZND1.imaginary.extend10+20i}}
	\caption{Solution $X_{\mathrm{R}}(\tau)$ computed by Con-CZND1 \cite{He2024ZeroingND} model and exact solution $X_{\mathrm{R}}^*(\tau)$ ($\gamma$ equals 10$+$20$\mathrm{i}$).}
	\label{fig.Con-CZND1.extend10+20i}
\end{figure*}

\begin{figure*}[!t]\centering
	\subfloat[]{\includegraphics[width=1\columnwidth]{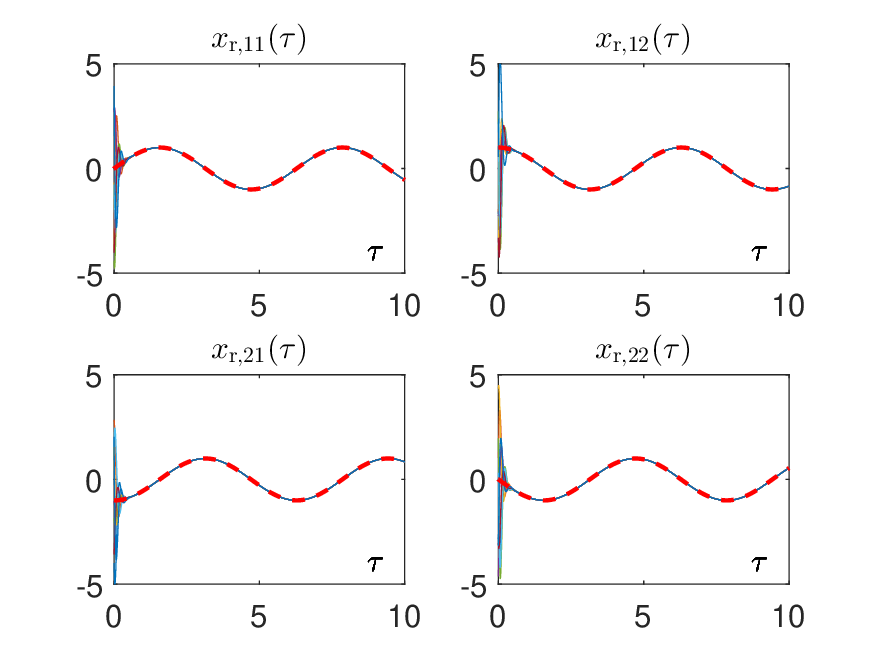}\label{fig.Con-CZND1_conj.real.extend10-20i}}
	\subfloat[]{\includegraphics[width=1\columnwidth]{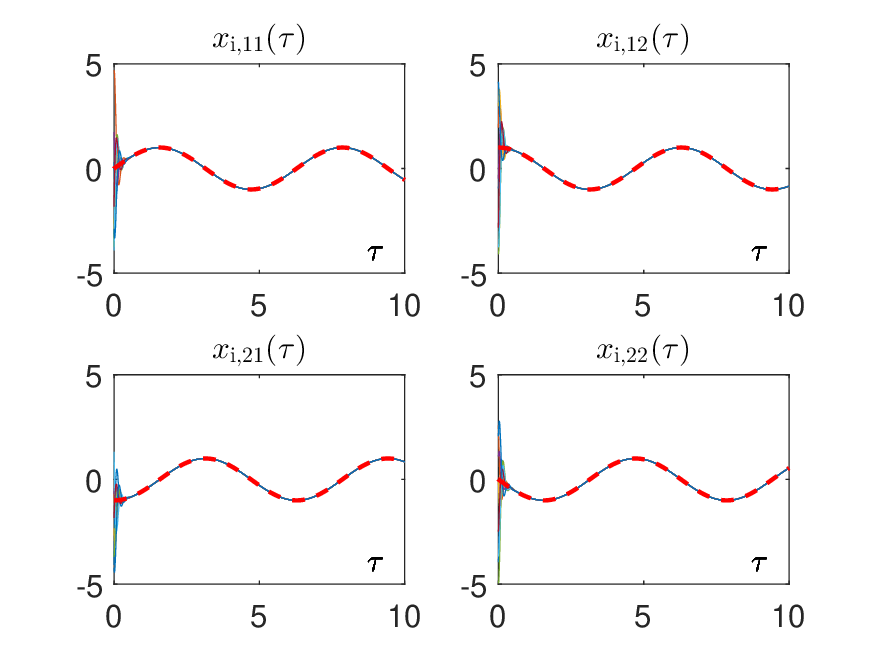}\label{fig.Con-CZND1_conj.imaginary.extend10-20i}}
	\caption{Solution $X_{\mathrm{R}}(\tau)$ computed by Con-CZND1\_conj \eqref{eq.solve.linearerrconcznd1.conj} model and exact solution $X_{\mathrm{R}}^*(\tau)$ ($\gamma$ equals 10$-$20$\mathrm{i}$).}
	\label{fig.Con-CZND1_conj.extend10-20i}
\end{figure*}

\begin{figure*}[!t]\centering
	\subfloat[]{\includegraphics[width=1\columnwidth]{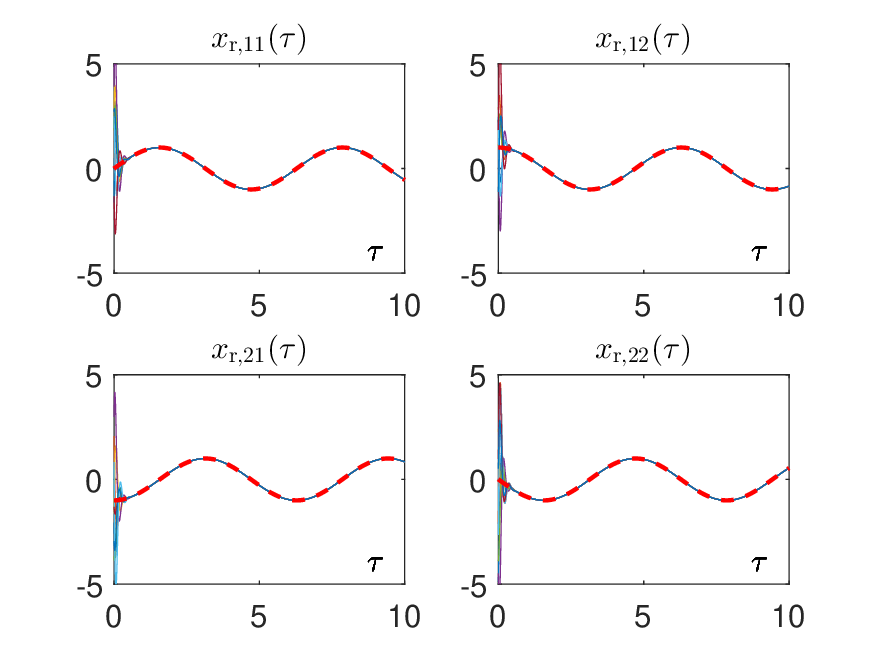}\label{fig.Con-CZND1.real.extend10-20i}}
	\subfloat[]{\includegraphics[width=1\columnwidth]{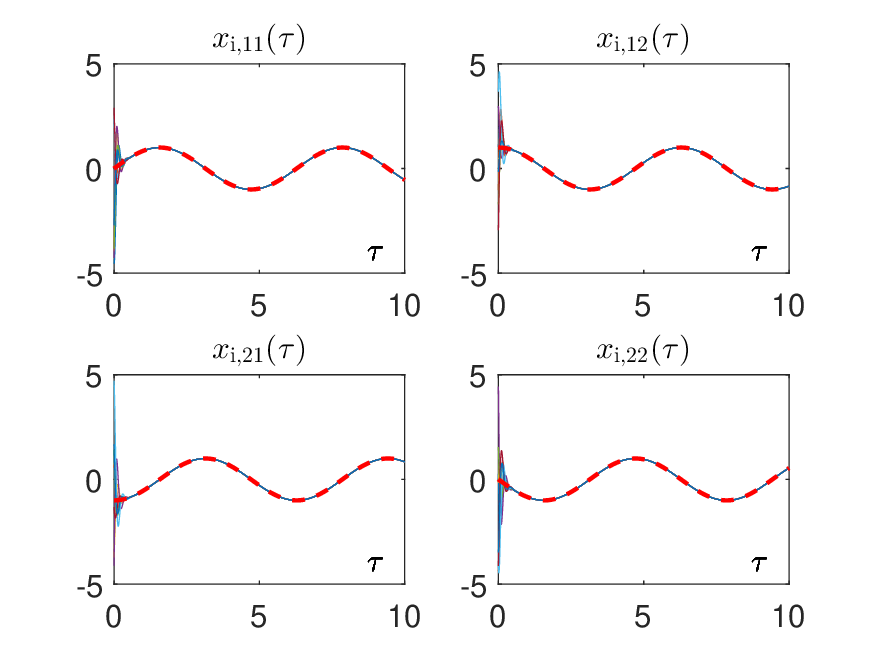}\label{fig.Con-CZND1.imaginary.extend10-20i}}
	\caption{Solution $X_{\mathrm{R}}(\tau)$ computed by Con-CZND1 \cite{He2024ZeroingND} model and exact solution $X_{\mathrm{R}}^*(\tau)$ ($\gamma$ equals 10$-$20$\mathrm{i}$).}
	\label{fig.Con-CZND1.extend10-20i}
\end{figure*}

\begin{figure*}[!t]\centering
	\subfloat[]{\includegraphics[width=1\columnwidth]{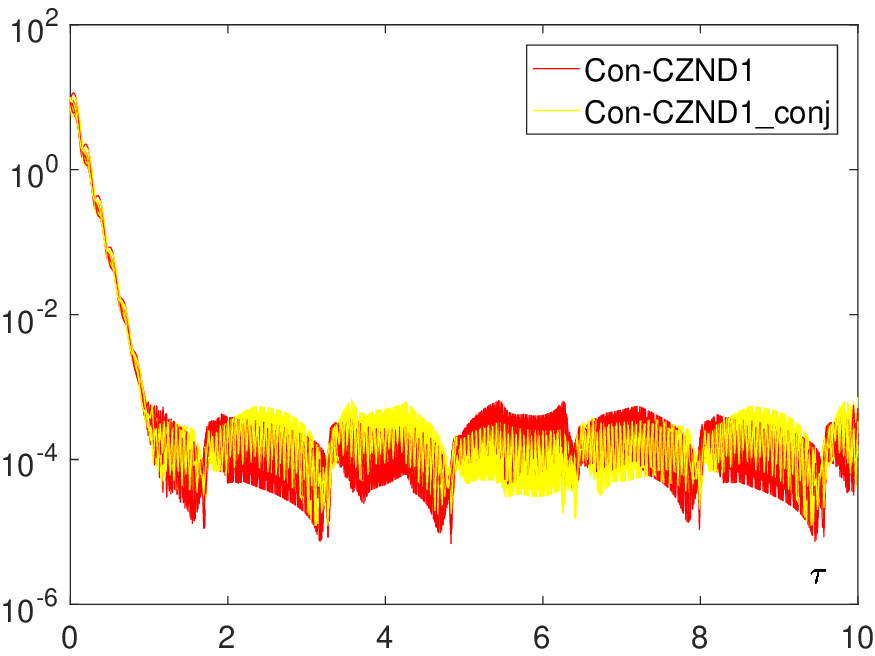}\label{fig.Con-CZND1_conj.vs.Con-CZND1.extend10+20i}}
	\subfloat[]{\includegraphics[width=1\columnwidth]{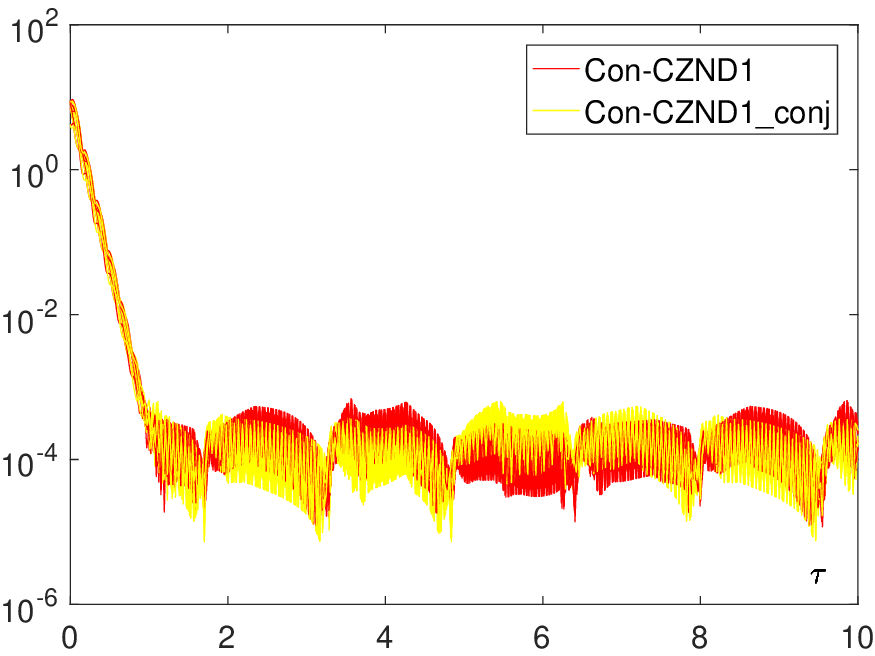}\label{fig.Con-CZND1_conj.vs.Con-CZND1.extend10-20i}}
	\caption{Logarithmic residual $\left \|X_{\mathrm{R}}(\tau)-X_{\mathrm{R}}^*(\tau)   \right \|_{\mathrm{F}}$ trajectories computed by Con-CZND1\_conj \eqref{eq.solve.linearerrconcznd1.conj} model vs. Con-CZND1 \cite{He2024ZeroingND} model. (a) $\gamma$ equals 10$+$20i. (b) $\gamma$ equals 10$-$20i.}
	\label{fig.Con-CZND1_conj.vs.Con-CZND1}
\end{figure*}

Since $E_{\mathrm{C2}}(\tau)\in\mathbb{C}^{m\times n}$ and $E_{\mathrm{C1}}(\tau)\in\mathbb{C}^{m\times n}$ are HN matrices. Then $\gamma$ is extended from $ \mathbb{R^+} $ to $\left (\mathbb{R^+} + \mathrm{i}\mathbb{R}\right )\in\mathbb{C}$, where 10$+$20$\mathrm{i}$ and 10$-$20$\mathrm{i}$ are selected. The running results of Con-CZND1\_conj \eqref{eq.solve.linearerrconcznd1.conj} model and Con-CZND1 model are shown in Figs. \ref{fig.Con-CZND1_conj.extend10+20i} through \ref{fig.Con-CZND1_conj.vs.Con-CZND1}.

As can be seen from Figs. \ref{fig.Con-CZND1_conj.extend10+20i} through \ref{fig.Con-CZND1.extend10-20i}, the imaginary part of $\gamma$ does not affect the results of using NHNSCAA, which is different from previous studies. Switching perspectives, the imaginary part of $\gamma$ is swallowed. Again, as seen in Fig. \ref{fig.Con-CZND1_conj.vs.Con-CZND1}, the imaginary part of $\gamma$ still makes Con-CZND1 model show some difference from its conjugate Con-CZND1\_conj \eqref{eq.solve.linearerrconcznd1.conj} model. More specific details of $\gamma$ runs can be found in {A\footnotesize{PPENDIX}}.
\subsection{Extension}
From the above study, it can be seen that HN matrices have a strong potential. However, the current theory of studying HN is not yet perfect. Therefore, in this paper, only complex matrices and TVSSCME \eqref{eq.TVSSCME} are skillfully used:
\begin{enumerate}
	\item{It is unavoidable to be transferred to RFTVLSLE.}
	\item{It is convenient to study the details of the space compressive approximation of NN.}
\end{enumerate}

To study the space compressive approximation error, there is an urgent and pressing need to develop the theory of HN, which is a key detail of NN black box study, and is also the current difficulty.

The following directions are recommended:
\begin{itemize}
	\item{Numerical linear algebra}
	\item{Differential geometry}
	\item{Qualitative and stability methods for ordinary differential equations}
	\item{Numerical solution of partial differential equations}
	\item{Nonlinear systems}
	\item{Switched systems}
	\item{Tensor optimizations}
	\item{Quantum mechanics}
	\item{Harmonic analysis}
	\item{Generalized functions theory}
\end{itemize}
\section*{Conclusion}
In this study, RFTVLSLE \eqref{eq.RFTVLSLE} corresponding to TVSSCME \eqref{eq.TVSSCME} is studied. More detailed theoretical explanations are given to highlight the two major ZND models, Con-CZND1 and Con-CZND2. Meanwhile, combining the above two major models, NHNSCAA is proposed from the perspective of space compressive approximation. On the basis of this approach, another Con-CZND1\_conj model is also constructed. With numerical experiments, the effectiveness of NHNSCAA is further supported. After, CCME will be further studied as the basis for solving NN black box.
\section*{Acknowledgments}
The authors thank Zhang et al. for their study of NN and Wu et al. for their study of CCME (of undivided importance). The authors join the studies of NN and CCME temporarily. It can be said that without their previous studies, the authors can not discover some keys.

Thanks again to them.

{\appendix[Con-CZND1 basic code]
\subsection{Main code start}
\begin{lstlisting}
function main()
close all
clear
clc
tspan=[0 10];
gamma= input('gamma = ');
random_length=10;%Length of random value.
Con_CZND1_Deal(gamma,tspan,random_length);
end
\end{lstlisting}\subsection{MatrixF}
\begin{lstlisting}
function output = MatrixF(t)
F=[6+sin(t),cos(t);cos(t),4+sin(t)]+1i*[cos(t),sin(t);sin(t),cos(t)];
output=F;
end
\end{lstlisting} 
\subsection{MatrixA} 
\begin{lstlisting}
function output = MatrixA(t)
A=[cos(t),sin(t);-sin(t),cos(t)]+1i*[sin(t),cos(t);cos(t),-sin(t)];
output=A;
end
\end{lstlisting} 
\subsection{MatrixC} 
\begin{lstlisting}
function output = MatrixC(t)
syms t;
C=[2*(cos(t))^2-2*cos(t)*sin(t)+6*sin(t),4*cos(t)+2*cos(t)*sin(t)-2*(cos(t))^2;-2*sin(2*t)-6*cos(t)+2,2*sin(2*t)-4*sin(t)-2]+1i*[2*(cos(t))^2+2*cos(t)*sin(t)+6*sin(t),4*cos(t)+2*cos(t)*sin(t)+2*(cos(t))^2;-2*sin(2*t)-6*cos(t)-2,-2*sin(2*t)-4*sin(t)-2];
output=C;
end
\end{lstlisting}
\subsection{ZNDTranLeftHandSide}
\begin{lstlisting}
function output = ZNDTranLeftHandSide(t,x)
if nargin<1,t=0;end;
[cm,cn]=size(MatrixC);
MatrixI_left=eye(cm);
MatrixI_right=eye(cn);
MatrixI_kron_A=kron(MatrixI_right,MatrixA(t));
MatrixF_T_kron_I=kron(conj(MatrixF(t))',MatrixI_left);
output=[real(MatrixF_T_kron_I)-real(MatrixI_kron_A),-imag(MatrixF_T_kron_I)-imag(MatrixI_kron_A);imag(MatrixF_T_kron_I)-imag(MatrixI_kron_A),real(MatrixF_T_kron_I)+real(MatrixI_kron_A)];
end
\end{lstlisting}
\subsection{AFMlinear}
\begin{lstlisting}
function Y = AFMlinear(X)
Y=X;
end
\end{lstlisting}
\subsection{ZNDTranRightHandSide}
\begin{lstlisting}
function output = ZNDTranRightHandSide(t,x,gamma)
if nargin==2,gamma=1;end;
syms u;
[cm,cn]=size(MatrixC);
F=diff(MatrixF(u));%Differentiate MatrixF.
A=diff(MatrixA(u));%Differentiate MatrixA.
C=diff(MatrixC(u));%Differentiate MatrixC.
u=t;
dotF=eval(F);%Execute the derivation statement.
dotA=eval(A);%Execute the derivation statement.
dotC=eval(C);%Execute the derivation statement.
MatrixC_eval=eval(MatrixC);%Execute.
for i=1:cm*cn,
x_dou(i)=complex(x(i),x(i+cm*cn));%Vectorized matrices with complex numbers.
end
MatrixX=reshape(conj(x_dou'),cm,cn);%Turn back to the original matrix containing the complex numbers, note that Matlab transpose the symbols encountered in the complex matrix into a conjugate transpose, so simply transpose the complex matrix then in the results of the conjugate once again.
disp(MatrixX);%Export MatrixX.
ZNDRightHandSide=dotC+dotA*conj(MatrixX)-(MatrixX)*dotF-gamma*AFMlinear(MatrixX*MatrixF(t)-MatrixA(t)*conj(MatrixX)-MatrixC_eval);%G complex matrix is operated.
ZNDRightHandSide_div=[real(ZNDRightHandSide),imag(ZNDRightHandSide)];%G complex matrix is extracted with the coefficients of the real and imaginary matrices.
v=reshape(ZNDRightHandSide_div,2*cm*cn,1);%The real and imaginary matrices are each vectorized and then combined into the form (Zr;Zi) and output.
output=v;
end
\end{lstlisting} 
\subsection{Con\_CZND1\_Deal}
\begin{lstlisting}
function Con_CZND1_Deal(gamma,tspan,random_length)
if nargin<1,gamma=1;end;
[cm,cn]=size(MatrixC);
options=odeset('Mass',@ZNDTranLeftHandSide,'MStateDep','none');%ODE setting.
for iter=1:8
x0=random_length*(rand(2*cm*cn,1)-0.5*ones(2*cm*cn,1));
[t,x]=ode45(@ZNDTranRightHandSide,tspan,x0,options,gamma);%ODE45 ordinary differential equation decoding.
%Plotting 2D diagrams with real and imaginary solutions shown separately.
f1=figure(1);
f1.Name='X real(Con-CZND1)';
f1.NumberTitle='off';
subplot(2,2,1);plot(t,x(:,1));title('$$x_{\mathrm{r},11}(\tau)$$','interpreter','latex');ylim([-random_length/2,random_length/2]);text('string','$$\tau$$', 'Units','normalized','HorizontalAlignment','center','Position',[0.90,0.10],'Interpreter','latex','FontSize',20,'FontWeight','Bold');  hold on %X11 real solution.
subplot(2,2,2);plot(t,x(:,1+cn));title('$$x_{\mathrm{r},12}(\tau)$$','interpreter','latex');ylim([-random_length/2,random_length/2]);text('string','$$\tau$$', 'Units','normalized','HorizontalAlignment','center','Position',[0.90,0.10],'Interpreter','latex','FontSize',20,'FontWeight','Bold');  hold on %X12 real solution.
subplot(2,2,3);plot(t,x(:,2));title('$$x_{\mathrm{r},21}(\tau)$$','interpreter','latex');ylim([-random_length/2,random_length/2]);text('string','$$\tau$$', 'Units','normalized','HorizontalAlignment','center','Position',[0.90,0.10],'Interpreter','latex','FontSize',20,'FontWeight','Bold');  hold on %X21 real solution.
subplot(2,2,4);plot(t,x(:,2+cn));title('$$x_{\mathrm{r},22}(\tau)$$','interpreter','latex');ylim([-random_length/2,random_length/2]);text('string','$$\tau$$', 'Units','normalized','HorizontalAlignment','center','Position',[0.90,0.10],'Interpreter','latex','FontSize',20,'FontWeight','Bold');  hold on %X22 real solution.
%---------------
f2=figure(2);
f2.Name='X imaginary(Con-CZND1)';
f2.NumberTitle='off';
subplot(2,2,1);plot(t,x(:,1+cm*cn));title('$$x_{\mathrm{i},11}(\tau)$$','interpreter','latex');ylim([-random_length/2,random_length/2]);text('string','$$\tau$$', 'Units','normalized','HorizontalAlignment','center','Position',[0.90,0.10],'Interpreter','latex','FontSize',20,'FontWeight','Bold');  hold on %X11 imaginary solution.
subplot(2,2,2);plot(t,x(:,1+cn+cm*cn));title('$$x_{\mathrm{i},12}(\tau)$$','interpreter','latex');ylim([-random_length/2,random_length/2]);text('string','$$\tau$$', 'Units','normalized','HorizontalAlignment','center','Position',[0.90,0.10],'Interpreter','latex','FontSize',20,'FontWeight','Bold');  hold on %X12 imaginary solution.
subplot(2,2,3);plot(t,x(:,2+cm*cn));title('$$x_{\mathrm{i},21}(\tau)$$','interpreter','latex');ylim([-random_length/2,random_length/2]);text('string','$$\tau$$', 'Units','normalized','HorizontalAlignment','center','Position',[0.90,0.10],'Interpreter','latex','FontSize',20,'FontWeight','Bold');  hold on %X21 imaginary solution.
subplot(2,2,4);plot(t,x(:,2+cn+cm*cn));title('$$x_{\mathrm{i},22}(\tau)$$','interpreter','latex');ylim([-random_length/2,random_length/2]);text('string','$$\tau$$', 'Units','normalized','HorizontalAlignment','center','Position',[0.90,0.10],'Interpreter','latex','FontSize',20,'FontWeight','Bold');  hold on %X22 imaginary solution.
end
%Precision solution.
%MatrixX is a 2 * 2 dimensional matrix.
x11_true=sin(t)+1i*sin(t); %X11 theoretical solution.
x12_true=cos(t)+1i*cos(t); %X12 theoretical solution.
x21_true=-cos(t)+1i*-cos(t);%X21 theoretical solution.
x22_true=-sin(t)+1i*-sin(t);%X22 theoretical solution.
%Plotting 2D diagrams with theoretical solution.
figure(1)
subplot(2,2,1);plot(t,real(x11_true),'r--','LineWidth',2);
subplot(2,2,2);plot(t,real(x12_true),'r--','LineWidth',2);
subplot(2,2,3);plot(t,real(x21_true),'r--','LineWidth',2);
subplot(2,2,4);plot(t,real(x22_true),'r--','LineWidth',2);
%---------------
figure(2)
subplot(2,2,1);plot(t,imag(x11_true),'r--','LineWidth',2);
subplot(2,2,2);plot(t,imag(x12_true),'r--','LineWidth',2);
subplot(2,2,3);plot(t,imag(x21_true),'r--','LineWidth',2);
subplot(2,2,4);plot(t,imag(x22_true),'r--','LineWidth',2);
\end{lstlisting}}

\bibliographystyle{IEEEtran}
\bibliography{bib-example}

\begin{IEEEbiography}[{\includegraphics[width=1in,height=1.25in,clip,keepaspectratio]{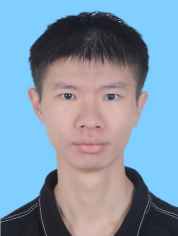}}]{Jiakuang He}
prefers the intersection of mathematics, chemistry, biology, medicine and computing.

Currently, he focuses on precision medicine and ODE neural networks based on complex conjugate matrix equations.
\end{IEEEbiography}

\begin{IEEEbiography}[{\includegraphics[width=1in,height=1.25in,clip,keepaspectratio]{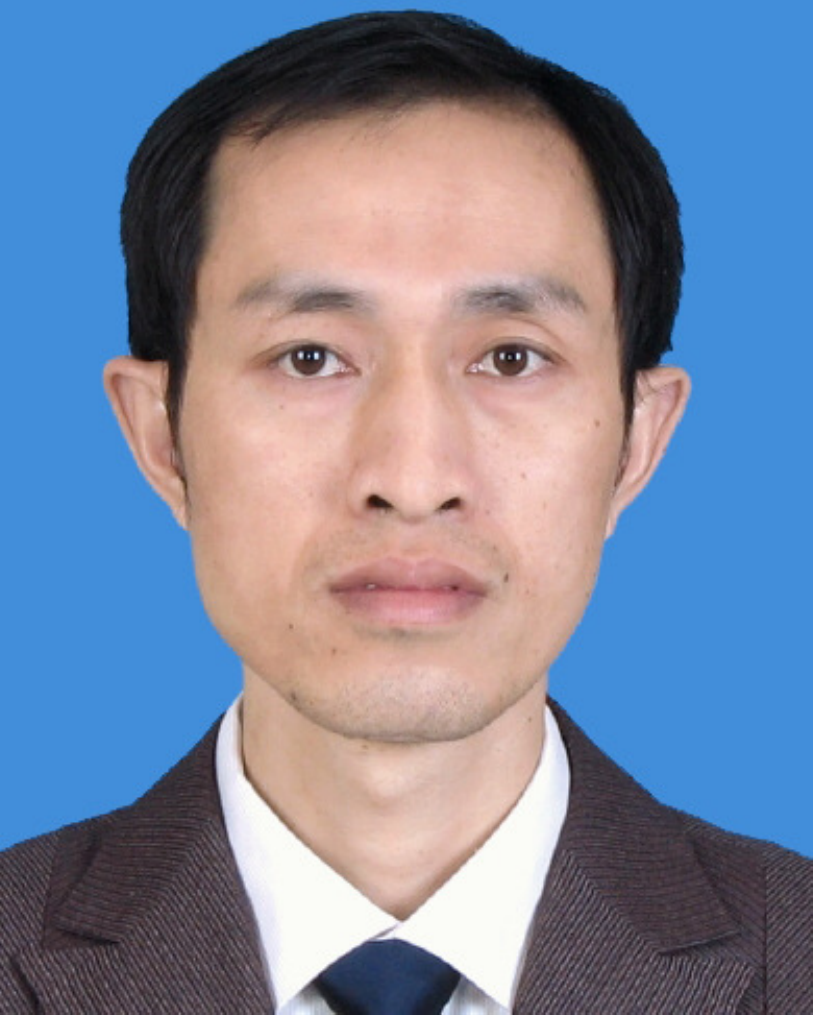}}]{Dongqing Wu}
	received the M.S. degree in computer graphics from the Institute of Industrial Design and Graphics, South China University of Technology, Guangzhou, China, in 2005, and the Ph.D. degree in mechanical engineering from the School of Electromechanical Engineering, Guangdong University of Technology, Guangzhou, China, in 2019. He is currently a professor with the School of Mathematics and Data Science, Zhongkai University of Agriculture and Engineering, Guangzhou, China. His current research interests include neural networks, robotics, and numerical analysis.
\end{IEEEbiography}
\vfill
\end{document}